%% file: main.tex
\documentclass[12pt]{amsart}

\usepackage{amsmath}
\usepackage{amsfonts}
\usepackage{amssymb}
\usepackage{amscd}
\usepackage{graphicx}
\usepackage[abbrev,alphabetic]{amsrefs}
\RequirePackage[dvipsnames,usenames]{color}
\usepackage{soul,xcolor}
\setstcolor{red}
\usepackage{stmaryrd}
\usepackage{mathtools}
\usepackage{booktabs}
\usepackage{multirow}
\usepackage{stmaryrd} 
\newtagform{tiny}{\tiny(}{)}

\usepackage{mathtools}
\usepackage{hyperref}
\usepackage[margin=1.25in]{geometry}

\usepackage{amsthm}
\usepackage{comment}
\usepackage[all,cmtip]{xy}
\usepackage{tikz-cd}
\usetikzlibrary{cd}
\usepackage{mathtools}
\usepackage[all]{xy}

\makeatletter
\def\@tocline#1#2#3#4#5#6#7{\relax
  \ifnum #1>\c@tocdepth 
  \else
    \par \addpenalty\@secpenalty\addvspace{#2}%
    \begingroup \hyphenpenalty\@M
    \@ifempty{#4}{%
      \@tempdima\csname r@tocindent\number#1\endcsname\relax
    }{%
      \@tempdima#4\relax
    }%
    \parindent\z@ \leftskip#3\relax \advance\leftskip\@tempdima\relax
    \rightskip\@pnumwidth plus4em \parfillskip-\@pnumwidth
    #5\leavevmode\hskip-\@tempdima
      \ifcase #1
       \or\or \hskip 1em \or \hskip 2em \else \hskip 3em \fi%
      #6\nobreak\relax
    \hfill\hbox to\@pnumwidth{\@tocpagenum{#7}}\par
    \nobreak
    \endgroup
  \fi}
\makeatother

\newsavebox{\pullback}
\sbox\pullback{%
\begin{tikzpicture}%
\draw (0,0) -- (1ex,0ex);%
\draw (1ex,0ex) -- (1ex,1ex);%
\end{tikzpicture}}

\newsavebox{\pullbackdl}
\sbox\pullbackdl{%
\begin{tikzpicture}%
\draw (-1ex,0ex) -- (0ex,0ex);%
\draw (0ex,-1ex) -- (0ex,0ex);%
\end{tikzpicture}}

\newsavebox{\pushoutdr}
\sbox\pushoutdr{%
\begin{tikzpicture}%
\draw (-1ex,-1ex) -- (-1ex,0ex);%
\draw (-1ex,0ex) -- (0ex,0ex);%
\end{tikzpicture}}

\def\phi{\varphi}
\def\epsilon{\varepsilon}
\def\tilde{\widetilde}

\def\mapsto{\longmapsto}

\newcommand{\rup}[1]{\lceil #1 \rceil}
\newcommand{\rdown}[1]{\lfloor #1 \rfloor}

\newcommand{\A}{\mathbb{A}}
\renewcommand{\P}{\mathbb{P}}
\newcommand{\Z}{\mathbb{Z}}
\newcommand{\Q}{\mathbb{Q}}

\newcommand{\R}{\mathbb{R}}
\newcommand{\F}{\mathbb{F}}

\newcommand{\cHom}{\mathcal{H}om}

\newcommand{\MO}{\mathcal{O}}
\newcommand{\sO}{\mathcal{O}}

\newcommand{\m}{\mathfrak{m}}

\newcommand{\univ}{\mathrm{univ}}
\newcommand{\red}{\mathrm{red}}

\newcommand{\Frac}{\mathrm{Frac}}
\newcommand{\Proj}{\mathrm{Proj}}
\renewcommand{\Im}{\mathrm{Im}}
\newcommand{\wt}{\widetilde}

\DeclareMathOperator{\gen}{gen}

\DeclareMathOperator{\Supp}{Supp}
\DeclareMathOperator{\Spec}{Spec}

\DeclareMathOperator{\Pic}{Pic}

\DeclareMathOperator{\Bs}{Bs}
\DeclareMathOperator{\Ex}{Ex}

\DeclareMathOperator{\Bl}{Bl}

\theoremstyle{plain}
\newtheorem{theorem}{Theorem}[section]
\newtheorem{thm}[theorem]{Theorem}

\newtheorem{prop}[theorem]{Proposition}

\newtheorem{lem}[theorem]{Lemma}

\newtheorem{cor}[theorem]{Corollary}

\newtheorem*{claim*}{Claim}
\newtheorem{step}{Step}

\theoremstyle{definition}

\newtheorem{dfn}[theorem]{Definition}
\newtheorem{definition}[theorem]{Definition}

\newtheorem{notation}[theorem]{Notation}
\newtheorem{nothing}[theorem]{}
\newtheorem*{setup*}{Setup}
\theoremstyle{plain}
\newtheorem{theo}{Theorem}

\theoremstyle{remark}
\newtheorem{rem}[theorem]{Remark}



\numberwithin{equation}{theorem}

\title[Weak quasi-F-splitting and del Pezzo varieties] 
{Weak quasi-$F$-splitting and del Pezzo varieties}

\author{Tatsuro Kawakami}
\address{Department of Mathematics, Graduate School of Science, Kyoto University, Kyoto 606-8502, Japan} 
\email{tatsurokawakami0@gmail.com}
\author{Hiromu Tanaka} 
\address{Graduate School of Mathematical Sciences, 
The University of Tokyo, 
3-8-1 Komaba, Meguro-ku, Tokyo 153-8914, JAPAN} 
\email{tanaka@ms.u-tokyo.ac.jp}

\begin{document}

\begin{abstract}
We show that smooth del Pezzo varieties in positive characteristic are quasi-$F$-split.
To this end, we introduce weak quasi-$F$-splitting and
we prove that general ladders of smooth del Pezzo varieties are normal. 
\end{abstract}

\maketitle

\setcounter{tocdepth}{2}

\tableofcontents
 \input section1.tex
 \input section2.tex

 \input section3.tex
 \input section4.tex
 \input section5.tex

 \input section6.tex

 \input main.bbl


\end{document}

%% file: section1.tex
\section{Introduction}
We work over an algebraically closed field $k$ of characteristic $p>0$. 
A variety $X$ is said to be \textit{$F$-split} if the Frobenius homomorphism  $\sO_X\to F_{*}\sO_X$ splits as an $\sO_X$-module homomorphism.
$F$-splitting is an important property of algebraic varieties in positive characteristic, as it assures Kodaira vanishing and the liftability to $W_2(k)$. 

Yobuko \cite{Yob19} recently introduced the notion of \textit{quasi-$F$-splitting}, which generalizes the notion of $F$-splitting (see also \cite{KTTWYY1} for fundamental results of quasi-$F$-splitting).
Quasi-$F$-split varieties still satisfy various useful properties that $F$-split varieties have, such as Kodaira vanishing 
and liftability to $W_2(k)$. 
Moreover, it is becoming clear that quasi-$F$-splitting is not as restrictive as $F$-splitting.
For example, all smooth del Pezzo surfaces are quasi-$F$-split although 
there exist smooth del Pezzo surfaces which are not $F$-split in low characteristics (\cite[Example 5.5]{Hara98}).

In this paper, we investigate quasi-$F$-splitting of \textit{del Pezzo varieties}, which form a special class of Fano varieties.
In the perspective of $\Delta$-genus (Definition \ref{def:delta genus}) introduced by Fujita \cite{Fujita75}, del Pezzo varieties can be seen as a natural generalization of del Pezzo surfaces to higher dimension.
\begin{dfn}\label{def:intro def}
    Let $(X,L)$ be a pair of a Gorenstein projective variety $X$ over an algebraically closed field $k$ and an ample Cartier divisor $L$.
    We say that $(X,L)$ is a \textit{del Pezzo variety} if the following hold. 
    \begin{enumerate}
        \item $-K_X\sim (n-1)L$, where $n\coloneqq \dim\,X$.
        \item $\Delta(X,L)=1$, where $\Delta(X, L)\coloneqq n +L^n-h^0(X, L)$. 
        \item $H^i(X,mL)=0$ for all $i\in\{1,\ldots,n-1\}$ and $m\in\Z$.
        \item $X$ is smooth along $\Bs\,|L|$.
    \end{enumerate}
    We call $L^{n}$ the \textit{degree} of $(X,L)$.
\end{dfn}
\begin{rem}
Assume that $X$ is smooth. 
If $n \leq 3$ or $k$ is of characteristic zero, 
then we only need (1) as the definition of del Pezzo varieties.
    See Remark \ref{rem:dPvar} for more details.
\end{rem}

Let $(X, L)$ be a del Pezzo variety. 
Then it is known that there exists a sequence, called \textit{ladder},  
\[
(X,L)
=:
(X_n,L_n) \supset (X_{n-1},L_{n-1}) \supset \cdots \supset (X_1,L_1)
\]
such that
\begin{enumerate}
    \item $(X_i,L_i)$ is a del Pezzo variety such that $\dim\,X_i=i$,
    \item $X_{i}\in |L_{i+1}|$ is a prime divisor on $X_{i+1}$, and 
    \item $L_{i}\sim L_{i+1}|_{X_i}$.
\end{enumerate}
When we take each $X_{i}\in |L_{i+1}|$ as a general member, 
such a sequence is called 
a \textit{general ladder}.
We refer to Theorem \ref{thm:ladder} for the existence of  ladders.

In this paper, 
we show that every smooth del Pezzo variety is quasi-$F$-split:

\begin{theo}[\textup{Theorem \ref{thm:normality of ladder} and Theorem \ref{mainthm:quasi-F-split}}]\label{introthm:rungs}
    Let $(X,L)$ be an $n$-dimensional smooth del Pezzo variety over an algebraically closed field $k$ of characteristic $p>0$.
    Let
    \[
   (X,L)=:
   (X_n,L_n) \supset (X_{n-1},L_{n-1}) \supset \cdots \supset (X_1,L_1)
    \] 
    be a general ladder.
    Then the following hold:
    \begin{enumerate}
        \item $X_i$ is normal for every integer $i\geq 1$.
        \item $X$ is quasi-$F$-split.
    \end{enumerate}
\end{theo}
\begin{rem}
    Smooth quasi-$F$-split varieties satisfy Kodaira vanishing (\cite[the proof of Theorem 4.1]{Yobuko2} or \cite[Theorem 3.15]{KTTWYY1})
    and there exist smooth Fano varieties that do not satisfy Kodaira vanishing (\cite{Tot19}).
    Therefore, smooth Fano varieties are not necessarily quasi-$F$-split. 
\end{rem}

We note that, in characteristic zero, every general ladder $X_i$ of a smooth del Pezzo variety $X$ is smooth. 
For the proof, we need Bertini's theorem for base point free linear systems, 
which is known to fail 
in positive characteristic.
Therefore, it is not trivial to show even the normality of $X_i$ in Theorem \ref{introthm:rungs}.

We also investigate the \textit{global $F$-regularity} of del Pezzo varieties.
Recall that global $F$-regularity is a stronger property than $F$-splitting in general, but it is equivalent to $F$-splitting for smooth Fano varieties (Remark \ref{rem:gfr}).

\begin{theo}[Section \ref{subsection:6-2}]\label{Intro:F-split}
    Let $(X,L)$ be an $n$-dimensional smooth del Pezzo variety 
    over an algebraically closed field $k$ of characteristic $p>0$ with $n \geq 2$. 
    Let \[
   (X,L)=: (X_n,L_n) \supset (X_{n-1},L_{n-1}) \supset \cdots \supset (X_1,L_1)
    \] be a general ladder.
    Then $X_i$ is globally $F$-regular for every integer $i\geq 2$ except for the following cases:
    \begin{enumerate}
        \item $L^n=3$ and $p=2$.
        \item $L^n=2$ and  $p \in \{2, 3\}$. 
        \item $L^n=1$ and  $p \in \{2, 3, 5\}$. 
    \end{enumerate}
\end{theo}

As mentioned before, due to the lack of Bertini's theorem, it is not easy to investigate singularities of $X_i$.
However,  
by the fact that globally $F$-regular varieties has only klt singularities, we can deduce from
Theorem \ref{Intro:F-split} that each $X_i$ 
has at worst canonical singularities when $p>5$ (Corollary \ref{cor:ladder is canonical}). 
For the proof of Theorem \ref{Intro:F-split}, we use global $F$-regularity of del Pezzo surfaces with canonical singularities, which is established in \cite{Kawakami-Tanaka(dPsurface)}.

\subsection{Sketch of proof of Theorem \ref{introthm:rungs}}

\subsubsection{Normality}
Fix a del Pezzo variety $(X,L)$ and a general ladder
\[
(X,L) =: (X_n,L_n) \supset (X_{n-1},L_{n-1}) \supset \cdots \supset (X_1,L_1).
\] 
The normality of each $X_i$ can be reduced to that of $X_1$. 
We here overview how to assure  
the smoothness of $X_1$ in the case $L^n=1$, which is the hardest case.
A key part is to establish 
the separability of the double cover  $\phi_{|2L|}\colon X\to \P(1,\dots,1,2)$ 
induced by the complete linear system $|2L|$ (see Lemma \ref{lem:degree 1}).
To show the separability, we generalize the earlier work by Megyesi \cite{Meg98} in dimension three to higher dimensional cases.

We will investigate smooth del Pezzo varieties of degree one in Section \ref{sec:dP1}.
In particular, we prove the separability of $\phi_{|2L|}$ in Subsection \ref{ss separability}.
We then prove the normality of general ladders (Theorem \ref{introthm:rungs}(1)) in Section \ref{sec:normality of ladder}.

\subsubsection{Quasi-$F$-splitting}
Take a general ladder  
\[
(X,L)=:(X_n,L_n) \supset (X_{n-1},L_{n-1}) \supset \cdots \supset (X_1,L_1). 
\]
By Theorem \ref{introthm:rungs}(1), each $X_i$ is normal. 
In particular, $X_1$ is smooth, and hence an elliptic curve. 
To clarify the idea, we start with the case when 
$X_1$ is $F$-split (i.e., an oridinary elliptic curve). 
In this case, 
it is easy to show that $X_2$ is $F$-split by using  an 
$F$-split inversion of adjunction. 
Repeating this procedure, 
we see that each $X_i$ (and hence $X$) is $F$-split.

Although $X_1$ is not necessarily $F$-split in general, 
$X_1$ is always quasi-$F$-split. 
The main strategy is to establish and apply a quasi-$F$-split inversion of adjunction.
A quasi-$F$-split inversion of adjunction 
for log Calabi-Yau `plt' pairs  is established in \cite[Theorem 4.6]{KTTWYY1}. 
However, this is not enough for our purpose, 
because we need to treat a pair which is not necessarily plt. 
To overcome this problem, 
we will introduce a new notion, which we call 
{\em weak quasi-$F$-splitting} (Section \ref{s-weak-QFS}). 
Then the strategy described in the previous paragraph is applicable after replacing 
$F$-splitting by weak quasi-$F$-splitting. 

As another technical obstruction, 
each $X_i$ is not necessarily smooth. 
In order to remedy this issue, 
we will use generic members instead of general members.  
Although generic members of base point free linear systems are 
always regular schemes (cf.\ \cite{Tan-Bertini}), 
they are defined over  an extension field $\kappa$ of $k$ which is usually an imperfect field (cf.\ Subsection \ref{ss generic Bertini}). 
For more details, see Section \ref{sec:quasi-F-split}. 
Weak quasi-$F$-splitting introduced in this article will 
play a crucial role for vanishing theorems and liftability 
of smooth Fano threefolds in positive characteristic 
\cite{KT-Fano3}.

\medskip
\noindent {\bf Acknowledgements.}
The authors thank Masaru Nagaoka, Teppei Takamatsu, Burt Totaro, and Shou Yoshikawa for valuable conversations related to the content of the paper.
Kawakami was supported by JSPS KAKENHI Grant number JP22KJ1771.
Tanaka was supported by JSPS KAKENHI Grant number JP22H01112 and JP23K03028. 

%% file: section2.tex
\section{Preliminaries}

\subsection{Notation and terminology}

In this subsection, we summarize notation and basic definitions used in this article. 
\begin{enumerate}
\item Throughout the paper, 
we work over an algebraically closed field $k$ of characteristic $p>0$ unless otherwise specified. 
\item We say that $X$ is a {\em variety} over a field $\kappa$ if 
    $X$ is an integral scheme 
    that is separated and of finite type over $\kappa$. 
    We say that $X$ is a {\em curve} (resp.~{\em surface}) 
    if $X$ is a variety of dimension one (resp.~two). 
\item For a variety $X$, 
we define the {\em function field} $K(X)$ of $X$ 
as the stalk $\sO_{X, \xi}$ at the generic point $\xi$ of $X$. 
\item We say that a finite surjective morphism $\phi\colon X\to Y$ of varieties is {\em separable} if the extension $K(X)/K(Y)$ of the function fields is separable.
\item Given a variety $X$ and a closed subscheme $Z$, 
we denote by $\Bl_Z X$ the blowup of $X$ along $Z$. 
\item Given a variety $X$, a closed subscheme $Z$, and a morphism $f\colon Y\to X$, 
we denote by $f^{-1}(Z)$ the scheme-theoretic inverse image of $Z$ by $f$, i.e., $f^{-1}(Z) := Y \times_X Z$. 
\item 
Given a projective variety $X$ and a Cartier divisor $L$ on $X$, 
$\Bs |L|$ denotes the scheme-theoretic base locus of an invertible sheaf $L$,  i.e., $\Bs |L|$ is the closed subscheme of $X$ 
whose ideal sheaf $I_{\Bs |L|}$ is given by 
\[
I_{\Bs |L|} \cdot \MO_X(L) 
= \Im( H^0(X, \MO_X(L)) \otimes_k \MO_X \to \MO_X(L)). 
\]
\item 
We say that $f: X \to Y$ is a {\em double cover} 
if $f$ is a finite surjective morphism of varieties such that 
the induced field extension $K(X) / K(Y)$ is of degree two. 
\item 
Let $X$ be an integral normal noetherian scheme. 
We say that a $\Q$-divisor $D$ is {\em simple normal crossing} 
if for every point $x \in \Supp D$, 
$X$ is regular at $x$ and 
there exist a regular system of parameters $x_1, ..., x_d$ of the regular local ring $\MO_{X, x}$ 
and $1 \leq r \leq d$ such that $\Supp D|_{\Spec \MO_{X, x}} = \Spec \MO_{X, x} / (x_1 \cdots x_r)$. Here $\Supp D$ denotes the reduced closed subscheme of $X$.
\end{enumerate}

\subsection{Bertini theorems for generic members}\label{ss generic Bertini}

In this  subsection, we gather some results on Bertini theorems for generic members.

\begin{definition}\label{d generic}
Let $K$ be a field and let $X$ be a projective variety over $K$, i.e., a projective integral scheme over $K$. 
Take an invertible sheaf $L$ on $X$ such that $H^0(X, L) \neq 0$. 
Set $X^{\gen}_L$ to be the {\em generic member} of 
$|L|$, i.e., for the family  $X^{\univ}_L$  parametrizing all the members in $|L|$, 
 we define $X^{\gen}_L$ by the following commutative diagram in which each square is cartesian (cf.\ \cite[Definition 5.6]{Tan-Bertini}): 
 \[
 \begin{tikzcd}
 X^{\gen}_L  \arrow[r] \arrow[d, hook] & X^{\univ}_L \arrow[d, hook]\\
 X \times_{K} K(\P(H^0(X, L))) \arrow[d] \arrow[r] & X \times_{K} \P(H^0(X, L))\arrow[d]\\
 \Spec K(\P(H^0(X, L)))\arrow[r] & \P(H^0(X, L)). 
 \end{tikzcd}
 \]
Note that the field extension $K \subset K(\P(H^0(X, L)))$ is a purely transcendental extension of finite degree. 
In \cite[Definition 5.6]{Tan-Bertini}, $X^{\gen}_L$ is written as $X^{\gen}_{L, 
H^0(X, L)}$. 
\end{definition}

\begin{nothing}\label{n generic repeated}
Fix an integer $m >0$. 
Let $k^{(0)}$ be a field and let $X^{(0)}$ be a  projective integral  scheme over $k^{(0)}$. 
Take an invertible sheaf $L^{(0)}$ on $X^{(0)}$ such that $H^0(X^{(0)}, L^{(0)}) \neq 0$. 
We shall construct  quadruples 
\[
\{(k^{(i)}, X^{(i)}, L^{(i)}, H^{(i)})\}_{1 \leq i \leq m}
\]
such that (1) and (2) below hold and 
\[
k^{(0)} \subset k^{(1)} \subset \cdots \subset k^{(m)}
\]
is a suitable sequence of field extensions. 
\begin{enumerate}
\item $X^{(i+1)} \coloneqq X^{(i)} \times_{k^{(i)}} k^{(i+1)} = X^{(0)} \times_{k^{(0)}} k^{(i+1)}$. 
\item $L^{(i+1)}$ is the pullback of $L^{(i)}$ (and hence of $L^{(0)}$). 
\end{enumerate}

Set $H^{(0)} \coloneqq 0$. 
Take an integer $0 \leq i <n$ and assume that 
the quadruple $(k^{(i)}, X^{(i)}, L^{(i)}, H^{(i)})$ has been already constructed. 
Then we set 
\begin{enumerate}
    \item[(3)] $k^{(i+1)} \coloneqq K(\P(H^0(X^{(i)}, L^{(i)})))$. 
\end{enumerate}
Then $X^{(i+1)}$ and $L^{(i+1)}$ are defined as in (1) and (2), respectively. 
\begin{enumerate}
    \item[(4)] Set $H^{(i+1)} \coloneqq (X^{(i)})_{L^{(i)}}^{\gen}$. 
    In particular, $H^{(i+1)}$ is a Cartier divisor on $X^{(i+1)}$ satisfying 
    $H^{(i+1)} \sim L^{(i+1)}$. 
\end{enumerate}
Repeating this procedure $m$ times, we obtain 
the quadruples $\{(k^{(i)}, X^{(i)}, L^{(i)}, H^{(i)})\}_{1 \leq i \leq m}$. 
Set 
\[
\kappa \coloneqq k^{(m)}, \qquad 
X_{\kappa} \coloneqq X^{(m)}, \qquad L_{\kappa}:=L^{(m)}. 
\]
Let $H_i$ be the effective Cartier divisor on $X_{\kappa} = X^{(m)}$ 
which is the pullback of $H^{(i)}$. 
By construction, the following hold. 
\begin{enumerate}
\renewcommand{\labelenumi}{(\roman{enumi})}
\item $k^{(0)} \subset \kappa$ is a purely transcendental extension of finite degree. 
\item $X_{\kappa} = X^{(0)} \times_{k^{(0)}} \kappa$. 
\item $L_{\kappa}$ is the pullback of $L^{(0)}$. 
\item $H_1, \ldots, H_m\in |L_{\kappa}|$. 
\end{enumerate}
\end{nothing}

\begin{lem}\label{l generic ladder}
Let $X$ be a smooth projective variety over $k$ of $\dim\,X=n$. 
Take an ample invertible sheaf $L$ such that one of (A) and (B) holds. 
\begin{enumerate}
\item[(A)] $|L|$ is base point free. 
\item[(B)] $L^n =1$ and a scheme-theoretic equality $\Bs |L|=P$ holds for a closed point $P$ of $X$. 
\end{enumerate}
For $m\coloneqq n-1$, $k^{(0)} \coloneqq k$, $X^{(0)}\coloneqq X$, and $L^{(0)} \coloneqq L$, we use the same notation as in 
(\ref{n generic repeated}). 
Then the following hold. 
\begin{enumerate}
\item $H_1+ \cdots +H_{n-1}$ is a simple normal crossing divisor on $X_{\kappa}$. 
\item  $H_{j_1} \cap \cdots \cap H_{j_i}$ is an integral scheme if $1 \leq i \leq n-1$ and $1 \leq j_1< \cdots <j_i\leq n-1$. 
\item $H_1 \cap \cdots \cap H_{n-1}$ is one-dimensional.  
\end{enumerate}
\end{lem}

\begin{proof}
We first treat the case  (A). 
Let $\varphi \coloneqq \varphi_{|L|} \colon X \to \P^n_k$ be 
the morphism induced 
by the complete linear system $|L|$ which is base point free. 
We then have $X^{\gen}_H = X^{\gen}_{\varphi}$ by \cite[Remark 5.8]{Tan-Bertini} (for the definition of $X^{\gen}_{\varphi}$, 
see \cite[Deftinition 4.3]{Tan-Bertini}). 
By \cite[Theorem 4.9(4)(12)]{Tan22}, $H^{(1)}$ is a regular prime divisor.  
Let $H^{(1)}_2$ be the pullback of $H^{(1)}$ 
by the projection $X^{(2)} = X^{(1)} \times_{k^{(1)}} k^{(2)} \to X^{(1)}$, 
i.e., $H^{(1)}_2 \coloneqq H^{(1)} \times_{k^{(1)}} k^{(2)}$. 
Since $k^{(1)} \subset k^{(2)}$ is a purely transcendental extension, 
$H^{(1)}_2$ is a regular prime divisor on $X^{(2)}$. 
By \cite[Remark 4.5]{Tan22}, $H^{(2)}|_{H^{(1)}_2}$ coincides with the generic member $(H^{(1)})^{\gen}_{\psi}$ of 
$\psi \colon H^{(1)} \hookrightarrow X^{(1)} \xrightarrow{\varphi^{(1)}} \P^n_{k^{(1)}}$, 
where $\varphi^{(1)} \colon X^{(1)} \to \P^n_{k^{(1)}}$ denotes $\varphi \times_k k^{(1)}$. 
Therefore, $H^{(2)}|_{H^{(1)}_2} (= H^{(1)}_2 \cap H^{(2)})$ is a regular prime divisor on $H^{(1)}_2$. 
Hence $H^{(1)}_2 + H^{(2)}$ is a simple normal crossing divisor with 
$\dim H^{(1)}_2 \cap H^{(2)} = \dim X-2$. 
For the pullbacks  $H^{(1)}_3$ and $H^{(2)}_3$ of $H^{(1)}$ and $H^{(2)}$ on $X^{(3)}$, 
the same argument implies that $H^{(1)}_3 + H^{(2)}_3 + H^{(3)}$ is 
a simple normal crossing divisor on $X^{(3)}$ satisfying $\dim (H^{(1)}_3 \cap H^{(2)}_3 \cap H^{(3)}) = \dim X-3$. 
By repeating this procedure, we are done for the case (A). 
The case (B) follows from a similar  argument to the one of the case (A) by using Lemma \ref{l deg1 regularity} below. 
\end{proof}

\begin{lem}\label{l deg1 regularity}
Let $K$ be an infinite field and let $Y$ be a regular projective variety over $K$ of $\dim\,Y=n$.  
Let $L$ be an ample invertible sheaf on $Y$.
Assume that $L^{n}=1$, $|L|$ is not base point free, and 
$\Bs |L|$ is zero-dimensional. 
Then the following hold. 
\begin{enumerate}
\item A scheme-theoretic equality $\Bs |L| = P$  holds for a $K$-rational point $P$. 
\item $Y$ is smooth at $P$. 
\item If $H\in |L|$ is an effective Cartier divisor, 
then $P \in H$ and $H$ is smooth at $P$. 
\item The generic member $Y^{\gen}_L$ of $|L|$ is regular. 
\end{enumerate}
\end{lem}

\begin{proof}
If $H_1$ is an arbitrary member of $|L|$ and 
$H_2,\ldots, H_{\ell} \in |L|$ are general members, 
then $H_1 \cap \cdots \cap H_{n}$ is zero-dimensional because so is $\Bs |L|$.  
Since $H_1 \cdot H_2 \cdots H_{n} = L^{n} = 1$, 
we get a scheme-theoretic equality $H_1 \cap \cdots \cap H_n = \Bs |L|=P$ for some $K$-rational point $P$.  
Thus (1) holds. 
Hence an arbitrary member $H(:=H_1)$ of $|L|$ is regular at $P$. 
Then each of $Y$ and $H$ is smooth at $P$ 
by \cite[Proposition 2.13]{Tan21i}. 
Thus (2) and (3) hold. 
By \cite[Theorem 4.9(4), Remark 5.8, Proposition 5.10]{Tan-Bertini}, 
$Y^{\gen}_L$ is regular outside the $K'$-rational point $P \times_K K'$, 
where $K' \coloneqq K( \P(H^0(Y, \MO_Y(H))))$.  
Hence (4) holds. 
\end{proof}

\begin{prop}\label{p generic vs general}    
Fix an integer $m>0$. 
Let $k$ be a field and let $X$ be a projective integral scheme over $k$. 
Take an invertible sheaf $L$ on $X$ such that $H^0(X, L) \neq 0$. 
For $k^{(0)} := k, X^{(0)} :=X$, and $L^{(0)} :=L$, we use the same notation as in (\ref{n generic repeated}). 
Set $P := \P(H^0(X, L))$. Let
\[
u : X^{\univ}_L \overset{j}{\hookrightarrow} X \times_k P \xrightarrow{{\rm pr}_2} P 
\]
be the universal family parametrizing the members of $|L|$. 
For $1 \leq i\leq m$ and the $m$-fold direct product $P^m := P \times_k \cdots \times_k P$, 
we define a closed subscheme $(X^{\univ}_L)_i$ of 
$X \times_k P^m$ by 
\[
(X^{\univ}_L)_i := \{ (x, z_1, ..., z_m) \in X \times P^m \,|\, 
(x, z_i) \in X^{\univ}_L, z_{i'} \in P\,\text{ for }\, i' \neq i\}. 
\]
Then the following hold. 
\begin{enumerate}
\item 
Take $D_1, ..., D_m \in |L|$. 
If  $[D_1, ..., D_m] \in P^m(k)$ denotes the $k$-rational point corresponding to the $m$-tuple $(D_1, ..., D_m)$, 
then $v^{-1}([D_1, ..., D_m]) \simeq D_1 \cap \cdots \cap D_m$. 
\item
$H_1 \cap \cdots \cap H_m$ is isomorphic to the generic fiber of 
\[
v : \bigcap_{i=1}^m (X^{\univ}_L)_i \hookrightarrow 
X \times_k P^m \xrightarrow{{\rm pr}} P^m. 
\]
\end{enumerate}
\end{prop}

\begin{proof}
(1) The assertion follows from construction.

(2) 
For simplicity, we only prove the case when $m=2$. 
Fix a $k$-linear basis: 
$H^0(X, L) = k\varphi_0 \oplus k\varphi_1 \oplus \cdots \oplus k \varphi_n$. 

We first compute the generic fiber of $v$. 
In order to distinguish the direct product factors of $P^2$, 
we set $P_1 := P$ and $P_2 :=P$, so that $P^2 = P \times_k P = P_1 \times_k P_2$. 
For $P_1$ and $P_2$, we introduce different coordinates: 
\[
P_1 =\P(H^0(X, L)) \xrightarrow{\simeq, \theta_1} \P^n_{[s_0: \cdots :s_n]} \qquad \text{and}\qquad
P_2 = \P(H^0(X, L)) \xrightarrow{\simeq, \theta_2} \P^n_{[t_0: \cdots :t_n]}. 
\]
We get 
\[
X^{\univ}_L \xrightarrow{\simeq, {\rm id} \times \theta_1} \{ s_0\varphi_0 + s_1 \varphi_1+ \cdots +s_n \varphi_n =0\} \subset X \times_k \P^n_{[s_0: \cdots :s_n]}.  
\]
\[
X^{\univ}_L \xrightarrow{\simeq, {\rm id} \times \theta_2} \{ t_0\varphi_0 + t_1 \varphi_1+ \cdots +t_n \varphi_n =0\} \subset X \times_k \P^n_{[t_0: \cdots :t_n]}.  
\]
Then 
\[
(X^{\univ}_L)_1 \cap (X^{\univ}_L)_2 
\xrightarrow{\simeq, {\rm id} \times \theta_1 \times \theta_2 }
\]
\[
\{ s_0\varphi_0 + s_1 \varphi_1+ \cdots +s_n \varphi_n =0\} \cap 
 \{ t_0\varphi_0 + t_1 \varphi_1+ \cdots +t_n \varphi_n =0\}
 \subset X \times_k \P^n_{[s_0: \cdots :s_n]} \times_k \P^n_{[t_0: \cdots :t_n]}. 
\]
For $\wt{s}_i := s_i/s_0$ and $\wt{t}_i := t_i/t_0$, 
the generic fiber $F$ of $v$ is given by 
\begin{equation}\label{e1 generic vs general} 
F \simeq 
\{ \varphi_0 + \wt{s}_1 \varphi_1+ \cdots +\wt{s}_n \varphi_n =0\} \cap 
 \{ \varphi_0 + \wt{t}_1 \varphi_1+ \cdots +\wt{t}_n \varphi_n =0\}
 \subset X \times_k K, 
\end{equation}
where $K := \Frac(\P^n_{[s_0: \cdots :s_n]} \times_k \P^n_{[t_0: \cdots :t_n]}) = k(\wt{s}_1, ..., \wt{s}_n, \wt{t}_1, ..., \wt{t}_n)$.

It is enough to prove that 
$H_1 \cap H_2$ is isomorphic to (\ref{e1 generic vs general}). 
Since $H^{(1)} := X^{\gen}_L$ is the generic fiber of $X^{\univ}_L \to P$, 
it holds that  
\[
H^{(1)} \simeq \{ \varphi_0 + \wt{s}_1 \varphi_1+ \cdots +\wt{s}_n \varphi_n =0\} \subset X \times_k k(\wt{s}_1, ..., \wt{s}_n).  
\]
Similarly, $H^{(2)}$ is given as follows: 
\[
H^{(2)} \simeq \{ \varphi_0 + \wt{t}_1 \varphi_1 + \cdots +\wt{t}_n \varphi_n =0\} \subset (X \times_k k(\wt{s}_1, ..., \wt{s}_n)) \times_{k(\wt s_1, ..., \wt s_n)} k(\wt s_1,..., \wt s_n, \wt t_1, ..., \wt t_n). 
\]
Then $\kappa = k(\wt s_1,..., \wt s_n, \wt t_1, ..., \wt t_n)$ and 
\[
H_1 \cap H_2 = 
 \{ \varphi_0 + \wt s_1 \varphi_1+ \cdots +\wt s_n \varphi_n =0\} 
 \cap  \{ \varphi_0 + \wt t_1 \varphi_1+ \cdots +\wt t_n \varphi_n =0\} \subset X \times_k \kappa. 
\]
Hence $H_1 \cap H_2$ is isomorphic to (\ref{e1 generic vs general}), as required. 
\end{proof}

\subsection{Del Pezzo varieties}
In this subsection, we gather basic facts on del Pezzo varieties.
In particular, we recall 
explicit descriptions of del Pezzo varieties of degree $d\geq 2$.

\begin{dfn}\label{def:delta genus}
    Let $(X,L)$ be a pair of a projective variety $X$ and an ample invertible sheaf $L$. 
    Set $n\coloneqq \dim\,X$. 
    We define 
    the \textit{$\Delta$-genus} $\Delta(X,L)$ of $(X,L)$ by
    \[
    \Delta(X,L)\coloneqq n+L^n-\dim H^0(X,L).
    \]
We say that a sequence 
\[
(X, L) =: (X_n, L_n) \supset (X_{n-1}, L_{n-1}) \supset \cdots 
\supset (X_1, L_1)
\]
is a {\em ladder} (of $(X, L)$) if 
\begin{enumerate}
\item $X_{i}$ is a member of $|L_{i+1}|$ which is an integral scheme, and  
\item $L_{i} := L_{i+1}|_{X_{i}}$
\end{enumerate}
for every $i \in \{ 1, 2, ..., n-1\}$. 
A {\em general ladder}  (of $(X, L)$) is a ladder of $(X, L)$ such that 
\begin{enumerate}
\item[(3)] $X_{i}$ is a general member of $|L_{i+1}|$ for every $i \in \{ 1, 2, ..., n-1\}$. 
\end{enumerate}
\end{dfn}

\begin{dfn}\label{d dP var}
    Let $(X,L)$ be a pair of a Gorenstein projective variety  $X$  and an ample invertible sheaf $L$.
    Set $n\coloneqq \dim\,X$. 
    We say that $(X,L)$ (or simply $X$) is a \textit{del Pezzo variety} if 
    \begin{enumerate}
        \item $-K_X\sim (n-1)L$.
        \item $\Delta(X,L)=1$.
        \item $H^i(X,mL)=0$ for all $i\in\{1,\ldots,n-1\}$ and $m\in\Z$.
        \item $X$ is smooth along $\Bs |L|$.
    \end{enumerate}
    We call  $L^{n}$ the \textit{degree} of $(X,L)$.
\end{dfn}
\begin{rem}\label{rem:dPvar}\,
   \begin{enumerate}
    \item Suppose that $X$ is smooth. 
    If $X$ satisfies Kodaira vanishing, 
    then (3) holds. 
    Moreover, (2) follows from (1) and (3) (an argument in the proof of \cite[Theorem 1.9]{Fujita(ch=0)} is valid in all characteristics). 
    Finally, (3) has been proven by Megyesi when $n=3$ 
    \cite[Lemma 5]{Meg98} (see also \cite[Corollary 3.6 and Corollary 3.7]{Kaw2} or \cite[Theorem 2.3]{FanoI} for Kodaira-type vanishing).
    \item 
    Our definition of 
    del Pezzo varieties
    (Definition \ref{def:delta genus}) 
    is more restrictive than the one by Fujita, while they coincide when $X$ is smooth. 
    Specifically, Fujita imposes the condition $g(X,L)=1$ instead of (4), {which is assured by (1), (2), and (4) (see the proof of Theorem \ref{thm:ladder}(2))}
          Although the opposite direction might not be true in general,
          the stronger condition (4) automatically holds 
          for all the situations appearing in this paper 
          {(cf.\ Theorem \ref{thm:ladder}(4)).}
   \end{enumerate}
   \end{rem}

\begin{thm}\label{thm:ladder}
Let $(X,L)$ be an $n$-dimensional  del Pezzo variety. 
Set $d\coloneqq L^n$.
Then the following hold. 
\begin{enumerate}
\item There exists a ladder of $(X, L)$. 
\item $g(X, L) =1$ (for the definition of $g(X, L)$, see \cite[(1.2)]{Fujita}). 
\item If $L^n \geq 2$, then $|L|$ is base point free. 
\item 
If a sequence
    \[
    (X,L)\coloneqq (X_n,L_n) \supset (X_{n-1},L_{n-1}) \supset \cdots \supset (X_1,L_1)
    \]
    is a ladder of $(X, L)$, then each 
    $(X_i,L_i)$ is a del Pezzo variety such that $\dim\,X_i=i$ and $L_i^{i}=d$. 
\end{enumerate}
\end{thm}

\begin{proof}
    Since $\Delta(X,L)=1$, it follows from 
    \cite[page 321, Theorem in Section 2]{Fujita} that 
    $\Bs |L|$ is a finite set.
    By \cite[Theorem 3.4]{Fujita}, 
    there exists a ladder 
    \[
    (X,L)=(X_n,L_n) \supset (X_{n-1},L_{n-1}) \supset \cdots \supset (X_1,L_1)
    \]
    of $(X, L)$. In particular, (1) holds. 
    
    Since $X_i$ is a Cartier divisor on $X_{i+1}$, the variety $X_i$ is Gorenstein for all $i\geq 1$.
    Then, by adjunction, we obtain $-K_{X_i}\sim (\dim\,X_{i}-1)L_i$ for all $i\geq 1$.
    By the definition of $L_i$, we have $L_i^i=L^n=d$.
    
    Let us show (2). 
    As proved in \cite[Section 1]{Fujita}, we have
    \[
    g(X,L)=g(X_{n-1},L_{n-1})=\cdots=g(X_1,L_1)=\dim H^1(X_1, \sO_{X_1}).
    \]
    Since $K_{X_1}\sim 0$ by adjunction, we have $\dim H^1(X, \sO_{X_1})=\dim H^0(X, \sO_{X_1})=1$ by Serre duality.
    Thus we obtain $g(X,L)=1$, and hence (2) holds. 
    Then the assertion (3) follows from \cite[Theorem 3.6(b)]{Fujita}.

Let us show (4). 
By \cite[Theorem 3.6 (a)]{Fujita}, we have 
    \[
    1=\Delta(X,L)=\Delta(X_{n-1},L_{n-1})=\cdots=\Delta(X_1,L_1). 
    \]
We now prove that $X_i$ is smooth along $\Bs |L_i|$. 
If $L^n \geq 2$, then $|L_i|$ is base point free by (3) and 
the surjectivity of $H^0(X, L) \to H^0(X_i, L_i)$ 
\cite[(1.5)]{Fujita}. 
It is enough to treat the case when $L^n=1$. 
In particular, $\Bs |L| =\Bs |L_i|$ is scheme-theoretically equal to 
reduced one point $P$, and hence $X_i$ is smooth along $P$ 
(Lemma \ref{l deg1 regularity}).

It suffices to prove  
    \[
    H^j(X_i,\sO_{X_i}(mL_{i}))=0
    \]
    for all $j\in\{1,\ldots,\dim\,X_i-1\}$ and $m\in\Z$.
    By descending induction on $i$, 
    this can be confirmed by the following short exact sequence: 
    \[
    0\to \sO_{X_{i+1}}((m-1)L_{j+1})\to \sO_{X_{i+1}}(mL_{i+1}) \to \sO_{X_{i}}(mL_{i})\to 0.
    \]
    
\end{proof}

\begin{prop}\label{prop:description}
    Let $(X,L)$ be an $n$-dimensional  smooth del Pezzo variety. 
    If $L^n\geq 3$, then $|L|$ is very ample.
\end{prop}
\begin{proof}
    The assertion follows from \cite[Corolary 5.3(c)]{Fujita}.
\end{proof}

\begin{thm}\label{thm:embedding(d=2)}
    Let $(X,L)$ be an $n$-dimensional smooth del Pezzo variety such that 
    $L^n=2$. Then the following hold:
    \begin{enumerate}
    \item $X$ can be embedded in a weighted projective space $\P(1,\ldots,1,2)_{[x_0:\cdots:x_{n}:y]}$ as a hypersurface
    \[
    X=\{y^2+f_2(x_0,\ldots,x_{n})y+f_4(x_0,\ldots,x_{n})=0\} \subset \P(1,\ldots,1,2)_{[x_0:\cdots:x_{n}:y]}
    \]
    of degree $4$, 
    where $f_2(x_0,\ldots,x_{n})$ (resp.~$f_4(x_0,\ldots,x_{n})$) is a homogeneous polynomial of degree 2 (resp.~4).
    \item $|L|$ is base point free and we have the following commutative diagram
    \[
    \xymatrix@C=60pt{
    X\ar[r]^{\phi_{|L|}}\ar@{^{(}->}[d] & \P^n \\
    \P(1,\ldots,1,2)\ar@{.>}[ru]_{\Pi}& , \\
    }
    \]
    where $\phi_{|L|}\colon X\to \P^n$ is 
    the morphism induced by $|L|$ and $\Pi$ is the rational map given by
    \[
    \Pi\colon \P(1,\ldots,1,2)_{[x_0:\cdots:x_{n}:y]}\dasharrow \P^n_{[x_0:\cdots:x_{n}]}\,\,[x_0 :\cdots : x_{n} :y] \mapsto [x_0 :\cdots : x_{n}].
    \]
    \item $\phi_{|L|}\colon X\to \P^n$ is a separable double cover.
    \end{enumerate}
\end{thm}
\begin{proof}
    (1)\,
    Essentially the same argument as in \cite[Theorem 8]{Meg98} shows that
    the section ring $R(X, L)=\bigoplus_{m\geq 0} H^0(X,mL)$ is generated by $x_0,\ldots,x_n\in H^0(X, L)$ and $y\in H^0(X, 2L)$, which prove the assertion.
    
    (2)\,
    By Theorem \ref{thm:ladder}(3), $|L|$ is base point free.  
    We note that 
    the morphism $\phi_{|L|}\colon X \to  \P^n$ sends $P\in X$ to $[x_0(P):\cdots:x_n(P)]$.
    Thus the diagram is commutative. 
    
    (3)\,
    Since $\deg(\phi_{|L|})=(\phi_{|L|}^{*}\sO_{\P^n}(1))^n=L^n=2$, the finite morphism $\phi_{|L|}$ is a double cover.
    We show that $\phi_{|L|}$ is separable. We may assume that $p=2$.
    If $\phi_{|L|}$ is not separable, it 
    is inseparable.
    However, an inseparable smooth double cover of $\P^n$ has to be a quadric hypersurface by \cite[Proposition 2.5]{Eke87}. 
    We note that a quadric hypersurface is not a del Pezzo variety as its Fano index is equal to $n$. 
\end{proof}

%% file: section3.tex
\section{Del Pezzo varieties of degree one}\label{sec:dP1}

The purpose of this section is to provide a detailed description of 
smooth del Pezzo varieties of degree one (Theorem \ref{thm:embedding(d=1)}). 

\subsection{Embedding to a weighted projective space}

We start by recalling some properties on Veronese varieties 
(Lemma \ref{lem:veronese}). 
Usually, a Veronese variety is defined as the image of $\P^m$ by $|\MO_{\P^m}(d)|$. 
In this article, we only use the case when $d=2$, and hence we adopt the following terminology.

\begin{dfn}\label{d-Vero}
We say that $V$ is the {\em $m$-dimensional Veronese variety} if 
$V$ is the image of the closed immersion, called the {\em Veronese embedding}, induced by 
the complete linear system $|\sO_{\P^m}(2)|$:
\[
\P^m \hookrightarrow \P^{ \frac{1}{2}(m+1)(m+2)-1}. 
\]
Note that  $h^0(\P^m, \sO_{\P^m}(2)) = \binom{m+2}{2} = \frac{1}{2} (m+1)(m+2)$. 
\end{dfn}

\begin{lem}\label{lem:veronese}
    Set  $Y := \P(1,\ldots,1,2)$  and $n:= \dim Y$. 
    Let $\beta\colon \tilde{Y}\to Y$ be the blow-up at $Q=[0:\cdots:0:1]$ and set $F\coloneqq \Ex(\beta)$. 
    Then the following hold.
    \begin{enumerate}
        \item 
        For $N \coloneqq n(n+1)/2$ and the closed immersion 
        \begin{eqnarray*}
        \iota : \P(1,\ldots,1,2) &\hookrightarrow& \P^N,\\
        {[ x_0: \cdots : x_{n-1} : y]} &\mapsto& {[x_0^2 : \cdots :x_{n-1}^2 :x_0x_1: \cdots :x_{n-2}x_{n-1}:y]}, 
        \end{eqnarray*}
        it holds that 
        \[
\iota(\P(1,\ldots,1,2)) = \widetilde V, 
        \]
        where $\widetilde V \subset \P^N$ denotes the projective cone 
        of the $(n-1)$-dimensional Veronese variety $V \subset \P^{N-1}$.
        \item $\tilde{Y}\simeq \P_{V}(\sO_V\oplus\sO_V(2))$, 
        where 
        we set $\MO_V(\ell) := (\theta^{-1})^*\MO_{\P^{n-1}}(\ell)$ via the isomorphism $\theta : \P^{n-1} \xrightarrow{\simeq} V$ induced by the Veronese embedding (Definition \ref{d-Vero}). 
        \item 
Let  $S_{+}$ (resp.~$S_{-}$) be the section 
of the projection $\pi : \widetilde{Y} = \P_V(\MO_V \oplus \MO_V(2)) \to V$
corresponding to $\sO_V\oplus \sO_V(2)\to \sO_V(2)$ 
(resp.~$\sO_V\oplus \sO_V(2)\to \sO_V$). 
Then $S_{+}\sim \sO_{\tilde{Y}}(1)$ and       
$S_{-}=F$.
        \item  $S_+  - S_- \sim 2\pi^*H_V$, where $H_V\in |\sO_V(1)|$ is a hyperplane.
        \item The scheme-theoretic equality $\beta^{-1}(Q) = F$ holds.
    \end{enumerate}
\end{lem}
\begin{proof}
(1)\, 
Take 
defining polynomials $f_1,\ldots, f_r$ of the $(n-1)$-dimensional  Veronese variety $V\subset \P^{\frac{1}{2}n(n+1)-1} = \P^{N-1}$ given as in 
\[
V= \Proj\,k[z_1,\ldots, z_N]/(f_1(z), \ldots, f_r(z)) \subset \P^{N-1}, 
\]
where $z$ denotes the multi-variable $(z_1, ..., z_N)$.
Recall that  the projective cone $\tilde{V} \subset \P^N$ over the Veronese variety $V\subset \P^{N-1}$ is defined by  
\[
\tilde{V} 
=\Proj\,k[z_1, \ldots, z_N, y]/(f_1(z), \ldots, f_r(z)) \subset \P^{N}. 
\]
By the construction of the Veronese variety, 
$[x_0^2 : \cdots :x_{n-1}^2 :x_0x_1: \cdots :x_{n-2}x_{n-1}:y]$ satisfies $f_1= \cdots =f_r=0$.
Thus the image of $\iota : \P(1,\ldots,1,2) \hookrightarrow \P^N$ is contained in $\widetilde{V}$. 
By $\dim \P(1, ..., 1, 2) =n= \dim \widetilde{V}$, we get the required equality  
$\iota(\P(1, ..., 1, 2)) = \widetilde V$.

(2)\, 
The assertion holds by applying \cite[Lemma 8.6]{FanoI} with $n\coloneqq \frac{1}{2}n(n+1)-1$, $s\coloneqq1$, and $\sO_X(1)\coloneqq \sO_V(2)$. 

(3), (4)\,
By the same argument as in \cite[Ch.\ V, 
Proposition 2.6]{Har77}, we get 
$\MO_{\widetilde Y}(1) \sim S_+$, $S_+|_{S_+} \sim 2\pi^*H_V|_{S_+}$, 
and $S_-|_{S_-} \sim -2\pi^*H_V|_{S_+}$. 
Since $S_+  - S_-$ is $\pi$-trivial, 
we can write $S_+  - S_-\sim \pi^*N$ for some $N\in \Pic(V)$. 
    We have $S_{+}|_{S_{+}}\sim \sO_{\tilde{Y}}(1)|_{S_{+}}\sim \sO_V(2)|_{S_{+}}$ and $S_{+}\cap S_{-}=\emptyset$, which implies  $N=S_{+}|_{S_{+}}=\sO_V(2)$. Thus (4) holds. 

It is enough to show $S_- =F$. 
It follows from $S_+|_{S_+} \sim 2\pi^*H_V|_{S_+}$ that $S_+$ is nef. 
Since $S_+  \sim  S_- +2\pi^*H_V$, we have that $S_+$ is big. 
Moreover, $S_+|_{S_-} \sim 0$, and hence $S_+$ is nef and big but not ample. 
Since $\rho(\tilde Y)=2$, it follows that $\tilde Y$ has exactly two extremal rays.
These rays correspond to $\pi\colon \tilde{Y} \to V$ and 
$\beta \colon \widetilde{Y} \to Y$. 
Since $S_{+}$ is nef, big, and not ample, it follows that $S_{+}$ corresponds to a birational morphism $\beta$, i.e., $S_+ \equiv a \beta^*H_Y$ for some rational number $a>0$ for $H_Y\in|\sO_Y(1)|$. 
Then both of $S_-$ and $F$ are prime divisors 
to which the restriction of $\beta^*H_Y$ 
is numerically trivial.
Since such a divisor has to be unique, we conclude $S_-=F$.

(5) By the same argument as in \cite[Lemma 8.6]{FanoI}, 
we have the following commutative diagram 
\[
\begin{tikzcd}
& \widetilde{Y}= \P_{V}( \MO_{V} \oplus \MO_{V}(2))=\Bl_Q \widetilde{V} 
\arrow[ld, "\pi"'] \arrow[rd, "\beta"] \arrow[dd, hook, "j_2"] &\\
V \arrow[dd, hook, "j_1"] && \widetilde{V}, \arrow[dd, hook, "j_3"]\\
&  \P_{\P^{N-1}}( \MO_{\P^{N-1}} \oplus \MO_{\P^{N-1}}(1)) =\Bl_Q \P^N \arrow[ld, "\pi'"'] \arrow[rd, "\beta'"] &\\
\P^{N-1} &&  \P^{N}, 
\end{tikzcd}
\]
where $\pi$ and $\pi'$ are the projections, $\beta$ and $\beta'$ are the blowups at $Q$, 
and the left square is cartesian. 
It is easy to see that the integral scheme $F' :=\Ex(\beta')$ is a section of $\pi'$ and we have a scheme-theretic equality $\beta'^{-1}(Q) = F'$. 
Then 
\[
\beta^{-1}(Q) = \beta^{-1}(j_3^{-1}(Q)) = j_2^{-1}(\beta'^{-1}(Q)) = j_2^{-1}(F') \overset{(\star)}{=} F, 
\]
where $(\star)$ is proven as follows. 
By $j_2^{-1}(F')_{\red} = F$, it is enough to show that $j_2^{-1}(F')$ 
is an integral scheme, 
which follows from the fact that 
$j_2^{-1}(F')$ is a section of $\pi$ (note that 
$F'$ is a section of $\pi'$ 
and the left square is cartesian.
\end{proof}

\begin{thm}\label{thm:embedding(d=1)}
    Let $(X,L)$ be an $n$-dimensional smooth del Pezzo variety such that 
    $L^n=1$. 
    Take the closed immersion: 
\begin{eqnarray*}
        \iota : \P(1,\ldots,1,2) &\hookrightarrow& \P^{\frac{n(n+1)}{2}}\\
        {[ x_0: \cdots : x_{n-1} : y]} &\mapsto& 
        {[x_0^2 : \cdots :x_{n-1}^2 :x_0x_1: \cdots :x_{n-2}x_{n-1}:y]}.  
        \end{eqnarray*}
    Then the following hold.
    \begin{enumerate}
    \item  $X$ can be embedded in a weighted projective space $\P(1,\ldots,1,2,3)_{[x_0:\cdots:x_{n}:y:z]}$ as a hypersurface
    \begin{multline*}
      \hspace{10mm}  X=\{z^2+f_3(x_0,\ldots,x_{n-1},y)z+f_6(x_0,\ldots,x_{n-1},y)=0\}\\ \subset \P(1,\ldots,1,2,3)_{[x_0:\cdots:x_{n-1}:y:z]}
    \end{multline*}
    of degree $6$, 
    where $f_3(x_0,\ldots,x_{n-1},y)$ (resp.~$f_6(x_0,\ldots,x_{n-1},y)$) is a 
    homogeneous polynomial of degree 3 (resp.~6).
    \item $|2L|$ is base point free, 
    $\varphi_{|2L|}(X) = \iota(\P(1, ..., 1, 2))$, and we have the following commutative diagram
    \[
    \xymatrix@C=60pt{
    X\ar[r]^{\phi}\ar@{^{(}->}[d] & \P(1,\ldots,1,2) \\
    \P(1,\ldots,1,2,3)\ar@{.>}[ru]_{\Pi}& , \\
    }
    \]
    where 
    \begin{itemize}
        \item $\phi\colon X\to \P(1,\ldots,1,2)$ is the morphism induced by $\varphi_{|2L|} : X \to \P^{\frac{n(n+1)}{2}}$ and $\varphi_{|2L|}(X) = \iota(\P(1, ..., 1, 2)) \simeq \P(1, ..., 1, 2)$, and 
        \item 
    $\Pi$ is the rational map given by 
    \begin{align*}
        \Pi\colon \P(1,\ldots,1,2,3)_{[x_0:\cdots:x_{n-1}:y:z]}&\dasharrow \P(1,\ldots,1,2)_{[x_0:\cdots:x_{n-1}:y]}\\
        [x_0 :\cdots : x_{n-1} :y:z] &\mapsto [x_0 :\cdots : x_{n-1}:y]
    \end{align*} 
    \end{itemize}
    \item $\phi$ is a double cover.
    \item $P:= \Bs |H|$ is a closed point.
           For the singular point $Q\coloneqq [0:\cdots:0:1] \in Y\coloneqq \P(1,\ldots,1,2)$, we have $\phi^{-1}(Q)_{\red} = P$ and $\phi^{-1}(Q) = \Spec\,\sO_{X, P}/\m_P^2$, 
           where $\m_P$ denotes the maximal ideal of $\MO_{X, P}$. 
    \item $\phi$ is separable.
    \end{enumerate}
\end{thm}

We now show (1)--(4) of Theorem \ref{thm:embedding(d=1)}. 
The proof of Theorem \ref{thm:embedding(d=1)}(5) will be given
in Subsection \ref{ss separability} (Corollary \ref{c-Chern-V1-p=2}).

\begin{proof}[Proof of (1)--(4) of Theorem \ref{thm:embedding(d=1)}]
    (1)\,
    Essentially the same argument as in \cite[Theorem 8]{Meg98} shows that
    the section ring $R(X, L)=\bigoplus_{m\geq 0} H^0(X,mL)$ is generated by $x_0,\ldots,x_{n-1}\in H^0(X, L)$, $y\in H^0(X, 2L)$, and $z\in H^0(X, 3L)$, which prove the assertion (1).

   (2)\,
   For $N := \frac{n(n-1)}{2}$, consider the following diagram: 
\[
\begin{tikzcd}
& \P^N\\
X \arrow[r, "\varphi"']\arrow[ru, "\varphi_{|2L|}"]  \arrow[d, hook] & \P(1, ..., 1, 2) \arrow[u, hook, "\iota"'] \\
 \P(1, ..., 1, 2,  3) \arrow[ru, dashrightarrow, "\Pi"'].
\end{tikzcd}
\]
It is enough to show the commutativity of this diagram. 
The lower triangle is commutative as follows: 
\[
\begin{tikzcd}
P \arrow[r, "\varphi"']  \arrow[d] & {[x_0(P): \cdots : x_{n-1}(P):y(P)]} \\
{[x_0(P): \cdots : x_{n-1}(P):y(P):z(P)].} \arrow[ru, "\Pi"']
\end{tikzcd}
\]
The upper triangle commutes by $H^0(X, 2L) = S^2H^0(X, L) \oplus ky$. 
Specifically, the mappings are given by 
\[
\begin{tikzcd}
& {[x_0(P)^2 : \cdots : x_{n-1}(P)^2 : x_{1}(P)x_2(P): \cdots : x_{n-2}(P)x_{n-1}(P) :y(P)]}\\
P \arrow[r, "\varphi"']\arrow[ru, "\varphi_{|2L|}"]  & {[x_0(P): \cdots : x_{n-1}(P):y(P)].}\arrow[u, "\iota"'] 
\end{tikzcd}
\]

    (3)\, 
    Recall that 
    $\iota(\P(1, ..., 1, 2)) = \widetilde{V} \subset \P^N$ 
    (Lemma \ref{lem:veronese}), 
    where $\widetilde V \subset \P^N$ is the projective cone over the $(n-1)$-dimensional Veronese variety $V \subset \P^{N-1}$. 
    We have $(\MO_{\P^N}(1)|_{\widetilde V})^n =  \deg \widetilde V = \deg V = 2^{\dim V} = 2^{n-1}$. 
Hence 
\[
2^n = (2L)^n = (\deg \phi)  \cdot (\MO_{\P^N}(1)|_{\widetilde V})^n = 
(\deg \phi) \cdot 2^{n-1}, 
\]
which implies $\deg \phi = 2$.

    (4)\,
    Since $x_0,\ldots,x_{n-1} \in |L|$ is a generator,
    the scheme-theoretic base scheme $\Bs |L|$ is defined by $\{x_0=\cdots=x_{n-1}=0\}$. 
    Since $L^n=1$, it follows that $P=\Bs |L|$ is a closed point, 
    and hence $x_0|_{\Spec \MO_{X, P}}, ..., x_{n-1}|_{\Spec \MO_{X, P}}$ is a generator of $\m_P$. 
   
Since 
    \[
Q' := \iota(Q) =     [0:\cdots :0:1]\in \mathrm{Im}(\phi)\subset \P^N,
    \]
    it follows that 
\[
\varphi^{-1}(Q) = \varphi^{-1}\iota^{-1}(Q') = \varphi_{|2L|}^{-1}( [0:\cdots :0:1])  
\]
\[
= \{x_0^2=\cdots=x_{n-1}^2 =x_0x_1= \cdots =x_{n-2}x_{n-1}=0\} 
= \Spec \MO_{X, P}/\m_P^2,
\]
and we conclude.
\end{proof}

\subsection{Separability}\label{ss separability}

This subsection is devoted to the proof of Theorem \ref{thm:embedding(d=1)}(5) (Corollary \ref{c-Chern-V1-p=2}).

\begin{prop}\label{p-Vero-cone}
We use the same notation as in Theorem \ref{thm:embedding(d=1)}. 
For the blowup $\alpha \colon \widetilde X\to X$ (resp.~$\beta \colon \widetilde Y \to Y$) at $P \in X$ (resp.~$Q \in Y$), 
there exists the following commutative diagram 
\[
\begin{tikzcd}
\widetilde{X} \arrow[r, "\alpha"] \arrow[d, "\widetilde \phi"] & X \arrow[d, "\phi"]\\
\widetilde{Y} \arrow[r, "\beta"]& Y
\end{tikzcd}
\]
such that $\widetilde \phi \colon \widetilde X \to \widetilde Y$ is a double cover.
\end{prop}
\begin{proof} 
Let $\overline{X}\coloneqq \mathrm{Bl}_{\phi^{-1}(Q)}X$.
Then by the universality of blow-up, there exists a unique morphism $\overline{X}\to \tilde{Y}$ which makes the following diagram commutative:
\[
\begin{tikzcd}
\overline{X} \arrow[r] \arrow[d, "\widetilde \phi"] & X \arrow[d, "\phi"]\\
\widetilde{Y} \arrow[r, "\beta"]& Y
\end{tikzcd}
\]
By Theorem \ref{thm:embedding(d=1)} (4), we have $\phi^{-1}(Q) = \Spec \MO_{X, P}/\m_P^2$.
Since $\m_P^2 \sO_{\tilde X} = (\m_P \sO_{\tilde X})^2 = \sO_{\tilde X}(-2E)$ is a Cartier divisor for $E \coloneqq \Ex(\alpha)$,
by the universality of blow-up again, $\alpha\colon \tilde{X}\to X$ factors through $\overline{X}\to X$.
Thus, we obtain the required commutative diagram.

 Set $F \coloneqq \Ex(\beta)$. 
Since $\widetilde{\phi}|_E \colon E (=\P^{n-1}) \to F(=\P^{n-1})$ is a surjective morphism of projective spaces, it have to be finite.
Therefore, $\widetilde \phi \colon \widetilde X \to \widetilde Y$ is a finite morphism, as it is quasi-finite and proper. 
Since $K(X) = K(\widetilde X)$ and $K(Y) = K(\widetilde Y)$, 
it follows from 
Theorem \ref{thm:embedding(d=1)}(3) that 
$\tilde \phi \colon \tilde X \to \tilde Y$ is a double cover.
\end{proof}

\begin{notation}\label{notation:deg=1}
We fix the following notation. 
\begin{enumerate}
    \item $(X, L)$ is a smooth del Pezzo variety such that $\dim\,X=n$ and $L^n =1$.
    \item $Y\coloneqq \P(1,\ldots,1,2)$. 
    \item $\phi\colon X\to Y$ is the double cover 
    induced by $\phi_{|2L|}$ (Theorem \ref{thm:embedding(d=1)}). 
    \item $P\coloneqq \Bs |L|$. 
    \item $Q\coloneqq [0:\cdots:0:1]\in Y$. 
    \item $\alpha\colon \tilde{X}\to X$ is the blow-up at $P$ with  $E\coloneqq \Ex(\alpha)$. 
    \item $\beta\colon \tilde{Y} \coloneqq \P_{V}(\sO_{V}\oplus \sO_V(2))\to Y$ is the blow-up at $Q$ with $F\coloneqq \Ex(\beta)$, where $ V\coloneqq  \P^{n-1}$.
    \item $\tilde{\phi}\colon \tilde{X}\to\tilde{Y}$ is the induced double cover (Proposition \ref{p-Vero-cone}). 
    \item $\tilde{B}$ is a Cartier divisor on $\tilde{Y}$ such that $K_{\tilde X} \sim \tilde{\phi}^*(K_{\tilde Y} +\tilde B)$ (\cite[Proposition 0.1.3]{CD89}) and $B\coloneqq \beta_{*}\tilde{B}$. 
    \item $\pi\colon \tilde{Y} = \P_{V}(\sO_{V}\oplus \sO_V(2))\to V$ is the projection (Lemma \ref{lem:veronese}). 
    \item $S_{+}$ (resp.~$S_{-}$) is the section corresponding to $\sO_V\oplus \sO_V(2)\to \sO_V(2)$ (resp.~$\sO_V\oplus \sO_V(2)\to \sO$). Note that $S_{+}\sim \sO_{\tilde{Y}}(1)$ and $S_{-} = F$ (Lemma \ref{lem:veronese}).
    \item Fix $H\in |L|$, $H_Y\in |\sO_Y(1)|$ and $H_V\in |\sO_V(1)|$, and set $\tilde{H}\coloneqq \pi^{*}H_V$.
          Recall that $H\sim \phi^{*}H_Y$ and $S_{+}-S_{-}\sim 2\tilde{H}$ (Lemma \ref{lem:veronese}).
    \item Set $S \coloneqq S_+ + S_-$ and $M \coloneqq 6\tilde{H}+4S_-\sim 2\tilde{H} + 2S$. 
\end{enumerate}
In particular, we have the following commutative diagram. 
\[
\begin{tikzcd}
S_{\pm} \arrow[r, hook] \arrow[d, "\widetilde{\phi}|_{S_{\pm}}", "\simeq"']&\widetilde{X} \arrow[r, "\alpha"] \arrow[d, "\widetilde \phi"] & X \arrow[d, "\phi"]\\
V & \P_{V}(\sO_{V}\oplus \sO_V(2)) =\widetilde{Y} \arrow[r, "\beta"]\arrow[l,"\pi"']& Y=\P(1, ..., 1, 2),
\end{tikzcd}
\]
\end{notation}

In Proposition \ref{p-tilde-L}, we will determine $\tilde B$ in Notation \ref{notation:deg=1} (9). 
To this end, we start with the following auxiliary result.

\begin{lem}\label{l-Vero-Weil}
We use Notation \ref{notation:deg=1}. 
Then the following hold. 
\begin{enumerate}
    \item $H_Y\sim\beta_{*}\tilde{H}$.
    \item $K_{\widetilde Y} \sim_{\Q} \beta^*K_Y + \frac{n-2}{2}F$. 
\end{enumerate}
\end{lem}

\begin{proof}
(1)\,Since $\pi\colon \tilde{Y} = \P_{V}(\sO_{V}\oplus \sO_V(2))\to V$ is a  $\P^1$-bundle and $S_-$ is a section of $\pi$, we have 
\[
\Pic(\tilde{Y})=\pi^{*}\Pic(V)\oplus \Z S_{-}=\Z\tilde{H}\oplus \Z S_{-}.
\]
Since $\Pic(\tilde{Y})=\mathrm{Cl}(\tilde{Y})\xrightarrow{\beta_*} \mathrm{Cl}(Y)=\Z H_Y$ is surjective, we have $H_Y\sim\beta_{*}\tilde{H}$.

(2)\,
We can write 
\[
K_{\widetilde Y} \sim_{\Q} \beta^*K_Y + aF
\]
for some $a\in\Q$.
By Lemma \ref{lem:veronese}(4), 
we get $\MO_F(-F) = \MO_{S_-}(-S_-) = \MO_{S_-}(2\pi^*H_V) = \MO_F(2)$.
Thus,
\[
\sO_F(-n) 
=K_F \sim (K_{\tilde Y} + F)|_F \sim \beta^*K_Y + (a+1)F|_F =(a+1)F|_F \sim \sO_F(-2(a+1)), 
\]
and we conclude.
\end{proof}

\begin{prop}\label{p-tilde-L}
We use Notation \ref{notation:deg=1}. 
Then we have
\[
\tilde{B} \sim 3\tilde{H} + 2F.
\]
\end{prop}

\begin{proof}
Recall that we have the following commutative diagrams, where the right one is obtained by restricting the left one: 
\[
\begin{tikzcd}
\widetilde{X} \arrow[r, "\alpha"] \arrow[d, "\widetilde \phi"] & X \arrow[d, "\phi"]\\
\widetilde{Y} \arrow[r, "\beta"]& Y
\end{tikzcd}
\hspace{20mm}
\begin{tikzcd}
E \arrow[r] \arrow[d] & P \arrow[d]\\
F\arrow[r]& Q. 
\end{tikzcd}
\]
We have an exact sequence 
\[
0 \to \sO_{\tilde{Y}} \to \tilde{\phi}_*\sO_{\tilde X} \to \sO_{\tilde{Y}} (-\tilde B) \to 0
\]
by \cite[Lemma A.1]{Kaw2}.
Then we have
\[
K_{\tilde X} \sim \tilde{\phi}^*(K_{\tilde Y} +\tilde{B})
\]
by \cite[Proposition 0.1.3]{CD89}.
Recall $B \coloneqq \beta_* \tilde{B}$ (Notation \ref{notation:deg=1} (9)).  
Note that $B$ is a Weil divisor, which is not necessarily Cartier. 
Applying the push-forward $\alpha_*$ for $K_{\tilde X} \sim \tilde{\phi}^*(K_{\tilde Y} +\tilde{B})$, 
we get
\[
K_X \sim \phi^*(K_Y + B). 
\]

Since $-K_X\sim (n-1)H$, $-K_Y\sim (n+2)H_Y$, and $H\sim\phi^{*}H_Y$, we have 
\[
\phi^{*}B\sim K_X-\phi^{*}K_Y\sim \phi^{*}(3H_Y), 
\]
and thus $B\sim 3H_Y$. 
It follows from Lemma \ref{l-Vero-Weil} that 
\[
\alpha^*K_X +(n-1)E \sim K_{\tilde X} \sim 
\tilde{\phi}^*(K_{\tilde Y}+\tilde B) 
\sim_{\Q} \tilde{\phi}^*( \beta^*K_Y + \frac{n-2}{2}F +\tilde B ).
\]
Note that 
\[
\tilde{\phi}^*F = \tilde{\phi}^{-1}F 
\overset{(\star)}{=}
\tilde{\phi}^{-1}\beta^{-1}Q = \alpha^{-1}\phi^{-1}Q = \alpha^{-1}(\Spec (\sO_{X, P}/\m_P^2)) = 2E, 
\]
where $(\star)$ follows from Lemma \ref{lem:veronese}(5). 
Therefore, 
\[
\alpha^*K_X  +E 
\sim_{\Q} \tilde{\phi}^*( \beta^*K_Y +\tilde{B} ) = 
\alpha^*\phi^*K_Y + \tilde{\phi}^*\tilde{B}. 
\]
Combining this with $K_X \sim \phi^*(K_Y+3 H_Y)$, we obtain
\[
\tilde{\phi}^* \tilde{B} \sim_{\Q} 
\alpha^*K_X  +E  -\alpha^*\phi^*K_Y 
\sim E + 3 \alpha^*\phi^*H_Y.
\]
Since $S_+  - S_- \sim 2\tilde{H}$ (Notation \ref{notation:deg=1}),  
we have $(\tilde{H} + \frac{1}{2}S_-)|_{S_-} \sim_{\Q} 0$, which implies $\beta^*H_Y \sim_{\Q} \tilde{H} + \frac{1}{2}S_- = \tilde{H} + \frac{1}{2}F$ by Lemma \ref{l-Vero-Weil} (1). 
Then 
\[
\tilde{\phi}^* \tilde{B} \sim_{\Q} E + 3 \alpha^*\phi^*H_Y \sim_{\Q}\tilde\phi^*(\frac{1}{2}F) +3 \tilde{\phi}^* (\tilde{H}+\frac{1}{2}F)  
= \tilde\phi^*( 3\tilde{H} + 2F).  
\]
Therefore, we get $\tilde B \sim_{\Q} 3\tilde{H} + 2F$, which implies 
$\tilde B \sim 3\tilde{H} + 2F$. 
This completes the proof.
\end{proof}

In what follows, we shall write $DE \coloneqq D \cdot E$ in order to ease the notation (e.g., $H^2S \coloneqq H \cdot H \cdot S$).

\begin{lem}\label{l-V1-p=2-Euler}
We use Notation \ref{notation:deg=1}. 
Then \begin{multline*}
    c_t(\Omega_{\tilde Y}^1 \otimes  \sO_{\tilde Y}(M) )(1 + (2\tilde{H}+2S)t)^2\\= (1 + (\tilde{H}+2S)t)^n (1+ (4\tilde{H}+3S)t+ (4\tilde{H}^2 + 6\tilde{H}S +2S^2)t^2). 
\end{multline*}
\end{lem}
\begin{proof}
 We have the following Euler sequence for toric varieties \cite[Theorem 8.1.6]{CLS11}: 
\begin{equation}\label{e1-V1-p=2-Euler}
0 \to \Omega^1_{X_{\Sigma}} \to \bigoplus_{\rho} \sO_{X_{\Sigma}}(-D_{\rho}) \to \Pic(X_{\Sigma}) \otimes_{\Z} \sO_{X_{\Sigma}} \to 0, 
\end{equation}
where $D_{\rho}$ in $\bigoplus_{\rho}$  runs over all the torus-invariant prime divisors. We apply this exact sequence by setting $\tilde{Y}= X_{\Sigma}$ 
(cf.~\cite[Example 7.3.5]{CLS11}). 
We have the following torus-invariant prime divisors. 
\begin{itemize}
\item The base $V=\P^{n-1}$ has the $n$ torus-invariant divisors $E_1,\ldots, E_n$, and hence its pullbacks $D_1, \ldots, D_n$ are torus-invariant. 
\item The sections $S_+$ and $S_-$ of $\pi \colon \tilde Y = \P_V(\sO_V \oplus \sO_V(2)) \to V$ are torus-invariant. 
\end{itemize}
It is easy to see that these are  all the torus-invariant prime divisors, 
because the complement of $S_+ \cup S_- \cup \bigcup D_i$ in $\tilde Y$ 
is a torus from which the torus action can be extended to $\tilde Y$. 
Recalling $\Pic(\tilde{Y})\simeq \Z^{\oplus 2}$, it follows from the Euler sequence (\ref{e1-V1-p=2-Euler}) $\otimes \sO_{\widetilde Y}(M)$ that
\[
 c_t(\Omega^1_{\tilde Y} \otimes \sO_{\tilde Y}(M) ) c_t( M)^2   =\prod_{D_{\rho} \in R} c_t(\sO(-D_{\rho}+M)), \qquad R \coloneqq \{D_1,\ldots, D_n, S_+, S_-\}. 
\]
We have $D_i \sim \pi^*H_V =\tilde{H}$ for every $i \in \{1, \ldots, n\}$, 
which implies 
\[
c_t(\Omega^1_{\tilde Y} \otimes  \sO_{\tilde Y}(M) ) (1+Mt)^2 = (1 + (M-\tilde{H})t)^n (1+ (M-S_+)t) (1+(M-S_-)t). 
\]
Since $M \sim 2\widetilde{H}+2S$, we obtain
\begin{itemize}
\item $1+Mt = 1+(2\tilde{H}+2S)t$, and 
\item $1 + (M-\tilde{H})t = 1 + (\tilde{H}+2S)t$. 
\end{itemize}
It follows from $S=S_++S_-$ and $S_+ \cap S_- = \emptyset$ that 
\begin{eqnarray*}
&&    (1+ (M-S_+)t) (1+(M-S_-)t)\\
&=&  1 + (2M - S) t + (M^2- MS)t^2\\
&=&  1 + (2(2\tilde{H}+2S) - S) t + ( (2\tilde{H}+2S)^2- (2\tilde{H}+2S)S)t^2\\
&=&  1 + (4\tilde{H}+3S) t + ( 4\tilde{H}^2 +6\tilde{H}S +2S^2)t^2. 
\end{eqnarray*}
Thus the assertion holds. 
\end{proof}

\begin{lem}\label{l-V1-p=2-HS-compute}
We use Notation \ref{notation:deg=1}. 
For every $k \in \{0, 1, 2,\ldots, n\}$, 
the following hold. 
\begin{enumerate}
\item $\tilde{H}^{n-k}  S^k = 2^{k-1} + (-2)^{k-1}$. 
\item $(\tilde{H}+2S)^{k}   \tilde{H}^{n-k} =\frac{5^{k}-(-3)^{k}}{2}$. 
\item $
(\tilde{H}+2S)^{k}   (-(2\tilde{H}+2S))^{n-k} =   \frac{5^{k}\cdot  (-6)^{n-k} - (-3)^{k} \cdot 2^{n-k}}{2}$. 
\end{enumerate}
\end{lem}

\begin{proof}
(1) 
We have $\tilde{H}^n =0$, which settles the case when $k=0$. 
We may assume $k \geq 1$. 
It follows from $S = S_+ +S_-$, $S_+ \cap S_- = \emptyset$, and $S_{\pm}|_{S_{\pm}} \simeq \pi^*(\pm 2 H_V)|_{S_{\pm}}$ that 
\begin{eqnarray*}
\tilde{H}^{n-k}   S^k &=&\tilde{H}^{n-k}  (S_+^k + S_-^k) \\
&=& (\tilde{H}|_{S_+})^{n-k}(S_+|_{S_+})^{k-1} + 
(\tilde{H}|_{S_-})^{n-k}  (S_-|_{S_-})^{k-1}\\
&=& H_V^{n-k}  (2H_V)^{k-1} + H_V^{n-k} (-2H_V)^{k-1}\\
&=& 2^{k-1} +( -2)^{k-1}. 
\end{eqnarray*}

(2) 
We have 
\[
(\tilde{H}+2S)^{k}   \tilde{H}^{n-k} = \sum_{\ell=0}^{k} \binom{k}{\ell} \tilde{H}^{k-\ell}   (2S)^{\ell} \tilde{H}^{n-k} \overset{{\rm (1)}}{=} 
\sum_{\ell=0}^{k} \binom{k}{\ell} 2^{\ell} (2^{\ell-1} + (-2)^{\ell-1}) 
\]
\[
=\frac{1}{2}\sum_{\ell=0}^{k} \binom{k}{\ell}  (4^{\ell} - (-4)^{\ell}) = \frac{1}{2}( (4+1)^{k} - (-4+1)^{k}) = \frac{5^{k}-(-3)^{k}}{2}. 
\]

(3) 
We have 
\begin{eqnarray*}
(\tilde{H}+2S)^{k}   (2\tilde{H}+2S)^{n-k} &=& 
(\tilde{H}+2S)^{k}   ((\tilde{H}+2S) +\tilde{H})^{n-k}\\
&=&  (\tilde{H}+2S)^k   \sum_{m=0}^{n-k} \binom{n-k}{m} (\tilde{H}+2S)^m \tilde{H}^{n-k -m}\\
&\overset{{\rm (2)}}{=}& \sum_{m=0}^{n-k} \binom{n-k}{m}\frac{5^{m+k}-(-3)^{m+k}}{2} \\
&=& 
\frac{5^{k}}{2}\sum_{m=0}^{n-k} \binom{n-k}{m}5^{m} -
\frac{(-3)^{k}}{2}\sum_{m=0}^{n-k} \binom{n-k}{m}(-3)^m\\
&=& \frac{5^{k}}{2} (5+1)^{n-k} -\frac{(-3)^{k}}{2} (-3+1)^{n-k} \\
&=& \frac{5^{k} \cdot 6^{n-k} - (-3)^{k} \cdot (-2)^{n-k}}{2}.     
\end{eqnarray*}
\end{proof}

\begin{lem}\label{l-V1-p=2-alpha-delta}
We use Notation \ref{notation:deg=1}. 
For $k \in \Z_{\geq 0}$, we set 
\begin{itemize}
\item $\alpha_k \coloneqq (k+1) (-(2\tilde{H}+2S))^k$, 
\item $\gamma_k\coloneqq \binom{n}{k} (\tilde{H}+2S)^k, \gamma_{-1} \coloneqq 0, \gamma_{-2} \coloneqq 0$, and 
\item $\delta_k \coloneqq \gamma_k + (4\tilde{H}+3S)\gamma_{k-1} + (4\tilde{H}^2+6\tilde{H}S+2S^2)\gamma_{k-2}$. 
\end{itemize}
Then 
\[
c_n(\Omega_{\tilde Y}^1 \otimes \sO_{\tilde Y}(M))  = \alpha_0 \delta_n + \alpha_1\delta_{n-1} + \cdots + \alpha_n \delta_0. 
\]
\end{lem}

\begin{proof}
Recall that  the following holds (Lemma \ref{l-V1-p=2-Euler}):
\[
c_t(\Omega_{\tilde Y}^1 \otimes M) (1 + (2\tilde{H}+2S)t)^2 = ( 1 + (\tilde{H}+2S)t)^n ( 1+ (4\tilde{H}+3S)t+ (4\tilde{H}^2+6\tilde{H}S+2S^2)t^2). 
\]
We have 
\[
( 1 + (\tilde{H}+2S)t)^n = \sum_{k=0}^n \binom{n}{k} (\tilde{H}+2S)^k t^k= \gamma_0 + \gamma_1 t + \gamma_2 t^2 + \cdots + \gamma_n t^n,
\]
which implies 
\begin{eqnarray*}
&&( 1 + (\tilde{H}+2S)t)^n ( 1+ (4\tilde{H}+3S)t+ (4\tilde{H}^2+6\tilde{H}S+2S^2)t^2)\\
&=& (\gamma_0+ \gamma_1 t + \gamma_2t^2 + \gamma_3t^3 + \cdots + \gamma_nt^n)  ( 1+ (4\tilde{H}+3)t+ (4\tilde{H}^2+6\tilde{H}S+2S^2)t^2)\\
&=& \gamma_0 + (\gamma_1 + (4\tilde{H}+3S)\gamma_0)t + (\gamma_2 + (4\tilde{H}+3S)\gamma_1 + (4\tilde{H}^2+6\tilde{H}S+2S^2)\gamma_0)t^2 + \cdots \\
&=& \sum_{k=0}^n (\gamma_k + (4\tilde{H}+3S)\gamma_{k-1} + (4\tilde{H}^2+6\tilde{H}S+2S^2)\gamma_{k-2}) t^k\\
&=& \sum_{k=0}^n \delta_k t^k. 
\end{eqnarray*}
Therefore, 
\[
c_t(\Omega_{\tilde Y}^1 \otimes \sO_{\tilde Y}(M)) 
= \left(\sum_{k=0}^n \delta_k t^k \right)
\cdot  \left( 1 -(2\tilde{H}+2S)t + (2\tilde{H}+2S)^2t^2 - (2\tilde{H}+2S)^3t^3+ \cdots \right)^2.  
\]

\begin{claim*}
$\left( 1 -(2\tilde{H}+2S)t + (2\tilde{H}+2S)^2t^2 - (2\tilde{H}+2S)^3t^3+ \cdots \right)^2 = {\displaystyle \sum_{k=0}^n \alpha_kt^k}$. 
\end{claim*}

\begin{proof}
Substituting $x= -(2\tilde{H}+2S)t$ to 
\[
(1+x+x^2 + \cdots)^2 =1+  2x + 3x^2 + 4x^3+ \cdots, 
\]
we obtain 
\[
\left( 1 -(2\tilde{H}+2S)t + (2\tilde{H}+2S)^2t^2 - \cdots \right)^2 = \sum_{k=0}^n (k+1) (-(2\tilde{H}+2S))^k t^k 
= \sum_{k=0}^n \alpha_k t^k. 
\]
This completes the proof of Claim. 
\end{proof}

We then get 
\[
c_t(\Omega_{\tilde Y}^1 \otimes \sO_{\tilde Y}(M)) 
= \left(\sum_{k=0}^n \alpha_k t^k \right)
\left(\sum_{k=0}^n \delta_k t^k \right), 
\]
which implies 
\begin{eqnarray*}
c_n(\Omega_{\tilde Y}^1 \otimes \sO_{\tilde Y}(M)) 
&=& \text{the coefficient of $t^n$ in }c_t(\Omega_{\tilde Y}^1 \otimes \sO_{\tilde Y}(M))\\
&=& \alpha_0\delta_n + \alpha_1\delta_{n-1} + \cdots +\alpha_n \delta_0. 
\end{eqnarray*}

\end{proof}

\begin{thm}\label{t-Chern-V1-p=2}
We use Notation \ref{notation:deg=1}. 
Then 
\[
c_n(\Omega_{\tilde Y}^1 \otimes \sO_{\tilde Y}(M)) = \frac{1}{3} ( 5^n-(-1)^n).
\]
\end{thm}

\begin{proof}
For an integer $m$ satisfying $m<0$ or $m>n$, 
we set $\alpha_m \coloneqq 0, \gamma_m \coloneqq 0$ and $\delta_m \coloneqq 0$. 
By Lemma \ref{l-V1-p=2-alpha-delta}, the following holds: 
\begin{eqnarray*}
&&c_n(\Omega_{\tilde Y}^1 \otimes \sO_{\tilde Y}(M))\\
&=& \sum_{k=0}^{\infty} \alpha_{k} \delta_{n-k}\\
&=& \sum_{k=0}^{\infty} (k+1)(-(2\tilde{H}+2S))^k ( \gamma_{n-k} + (4\tilde{H}+3S)\gamma_{n-k-1} +(4\tilde{H}^2+6\tilde{H}S+2S^2)\gamma_{n-k-2})\\
&=& \sum_{k=0}^{\infty} (k+1)(-(2\tilde{H}+2S))^k \gamma_{n-k}\\
&+& (4\tilde{H}+3S) \sum_{k=0}^{\infty} (k+1)(-(2\tilde{H}+2S))^k \gamma_{n-k-1}\\
&+&  (4\tilde{H}^2+6\tilde{H}S+2S^2) \sum_{k=0}^{\infty} (k+1)(-(2\tilde{H}+2S))^k \gamma_{n-k-2}\\
&=& \sum_{k=0}^{\infty} (k+1) \binom{n}{n-k} (\tilde{H}+2S)^{n-k}(-(2\tilde{H}+2S))^k\\
&+& (4\tilde{H}+3S) \sum_{k=0}^{\infty} (k+1) \binom{n}{n-k-1} (\tilde{H}+2S)^{n-k-1}(-(2\tilde{H}+2S))^k\\
&+&  (4\tilde{H}^2+6\tilde{H}S+2S^2)  \sum_{k=0}^{\infty} (k+1) \binom{n}{n-k-2} (\tilde{H}+2S)^{n-k-2}(-(2\tilde{H}+2S))^k.
\end{eqnarray*}

\medskip

\setcounter{step}{0}

\begin{step}\label{s1-Chern-V1-p=2}
We have 
\[
\sum_{k=0}^{\infty} (k+1) \binom{n}{n-k} (\tilde{H}+2S)^{n-k}(-(2\tilde{H}+2S))^k = 4n(-1)^n.
\]
\end{step}

\begin{proof}[Proof of Step \ref{s1-Chern-V1-p=2}]
Since $(1+z)^n = \sum_{k=0}^{\infty}  \binom{n}{k} z^k$, 
we have
\[ 
(1+z)^n + nz(1+z)^{n-1}= \frac{d}{dz}(z(1+z)^n) = \sum_{k=0}^{\infty}  (k+1) \binom{n}{k} z^k. 
\]
Set $ z\coloneqq y/x$ and take the multiplication with $x^n$: 
\[
(x+y)^n + ny(x+y)^{n-1} = \sum_{k=0}^{\infty}  (k+1) \binom{n}{n-k} x^{n-k}y^k. 
\]
By setting  $x \coloneqq \tilde{H}+2S$ and $y \coloneqq -(2\tilde{H}+2S)$, we obtain $x+y = -\tilde{H}$, and hence  
\[
 \sum_{k=0}^{\infty} (k+1) \binom{n}{n-k} (\tilde{H}+2S)^{n-k}(-(2\tilde{H}+2S))^k
 \]
 \[
 = (-\tilde{H})^n -n(2\tilde{H}+2S)(-\tilde{H})^{n-1} = 0 +(-1)^n2n S\tilde{H}^{n-1} = (-1)^n 4n,
\]  
where we used $\tilde{H}^n =0$ and $\tilde{H}^{n-1}S = 2$ (Lemma \ref{l-V1-p=2-HS-compute}). 
This completes the proof of Step \ref{s1-Chern-V1-p=2}. 
\end{proof}

\begin{step}\label{s2-Chern-V1-p=2}
We have 
\[
(4\tilde{H}+3S) \sum_{k=0}^{\infty} (k+1) \binom{n}{n-k-1} (\tilde{H}+2S)^{n-k-1}(-(2\tilde{H}+2S))^k = 6n(-1)^{n-1}.
\]
\end{step}

\begin{proof}[Proof of Step \ref{s2-Chern-V1-p=2}]
We have 
\begin{align*}
    \sum_{k=0}^{\infty} (k+1) \binom{n}{k+1} z^k  = 
\frac{d}{dz}\left(\sum_{k=0}^{\infty} \binom{n}{k+1} z^{k+1}\right)  =& 
\frac{d}{dz}\left(\sum_{\ell=0}^{\infty} \binom{n}{\ell} z^{\ell}\right) \\
=& \frac{d}{dz}( (1+z)^n) = n (1+z)^{n-1}. 
\end{align*}
Set  $z \coloneqq y/x$ and take the multiplication with $x^{n-1}$: 
\[
\sum_{k=0}^{\infty} (k+1) \binom{n}{n-k-1} x^{n-k-1}y^k = n(x+y)^{n-1}. 
\]
Therefore, 
\[
(4\tilde{H}+3S) \sum_{k=0}^{\infty} (k+1) \binom{n}{n-k-1} (\tilde{H}+2S)^{n-k-1}(-(2\tilde{H}+2S))^k
\]
\[
= (4\tilde{H}+3S) \times n (-\tilde{H})^{n-1} = (-1)^{n-1}3 n S\tilde{H}^{n-1} =  6n(-1)^{n-1}.  
\]
where we used $\tilde{H}^n=0$ and $\tilde{H}^{n-1}S = 2$ (Lemma \ref{l-V1-p=2-HS-compute}). 
This completes the proof of Step \ref{s2-Chern-V1-p=2}. 
  \end{proof}

\noindent 
What is remaining is to compute 
\[
(4\tilde{H}^2+6\tilde{H}S+2S^2) \sum_{k=0}^{\infty} (k+1) \binom{n}{n-k-2} (\tilde{H}+2S)^{n-k-2}(-(2\tilde{H}+2S))^k.  
\]

\begin{step}\label{s3-Chern-V1-p=2}
For $x \in \R \setminus \{0\}$, we have
\[
f(x) \coloneqq \sum_{k=0}^{\infty} (k+1) \binom{n}{n-k-2} x^{k} =  \frac{nx (x+1)^{n-1} -(x+1)^n +1}{x^2}.
\]
\end{step}

\begin{proof}[Proof of Step \ref{s3-Chern-V1-p=2}]
\[
f(x) = \sum_{k=0}^{\infty} (k+1) \binom{n}{n-k-2} x^{k}
= \left( \sum_{k=0}^{\infty} \binom{n}{k+2} x^{k+1}\right)'
= \left( \sum_{\ell=0}^{\infty} \binom{n}{\ell+1} x^{\ell}\right)'
\]
\[
= \left( \frac{\sum_{\ell=0}^{\infty} \binom{n}{\ell+1} x^{\ell+1}}{x}\right)'
= \left( \frac{-1+\sum_{m=0}^{\infty} \binom{n}{m} x^{m}}{x}\right)'
= \left( \frac{-1+(x+1)^n}{x}\right)'
\]
\[
= \frac{nx (x+1)^{n-1} -(x+1)^n +1}{x^2}. 
\]
This completes the proof of Step \ref{s3-Chern-V1-p=2}. 
\end{proof}

\begin{step}\label{s4-Chern-V1-p=2}
We have
\[
(4\tilde{H}^2+6\tilde{H}S +2S^2) \sum_{k=0}^{\infty} (k+1) \binom{n}{n-k-2} (\tilde{H}+2S)^{n-k-2}(-(2\tilde{H}+2S))^k 
\]
\[
=  \frac{1}{3}((-1)^n(6n-1)+5^n). 
\]
\end{step}

\begin{proof}[Proof of Step \ref{s4-Chern-V1-p=2}]
It holds that 
\begin{eqnarray*}
&&(4\tilde{H}^2+6\tilde{H}S +2S^2) \sum_{k=0}^{\infty} (k+1) \binom{n}{n-k-2} (\tilde{H}+2S)^{n-k-2}(-(2\tilde{H}+2S))^k\\
&=& 2(2\tilde{H}+S)(\tilde{H}+S) \sum_{k=0}^{\infty} (k+1) \binom{n}{n-k-2} (\tilde{H}+2S)^{n-k-2}(-(2\tilde{H}+2S))^k\\
&=& ((\tilde{H}+2S) - (3\tilde{H}+3S)) (-2\tilde{H}-2S) \sum_{k=0}^{\infty} (k+1) \binom{n}{n-k-2} (\tilde{H}+2S)^{n-k-2}(-(2\tilde{H}+2S))^k\\
&=&  {\rm (i)} + \frac{3}{2} \times {\rm (ii)}, 
\end{eqnarray*}
where 
\begin{eqnarray*}
{\rm (i)} &\coloneqq & \sum_{k=0}^{\infty} (k+1) \binom{n}{n-k-2} (\tilde{H}+2S)^{n-k-1}(-(2\tilde{H}+2S))^{k+1}\\
&=&\sum_{k=0}^{\infty} (k+1) \binom{n}{n-k-2} \frac{5^{n-k-1} \cdot (-6)^{k+1} -(-3)^{n-k-1}\cdot 2^{k+1}}{2}
\end{eqnarray*} 
\begin{eqnarray*}
{\rm (ii)}  
& \coloneqq& \sum_{k=0}^{\infty} (k+1) \binom{n}{n-k-2} (\tilde{H}+2S)^{n-k-2}(-(2\tilde{H}+2S))^{k+2}\\
&=&\sum_{k=0}^{\infty} (k+1) \binom{n}{n-k-2} \frac{5^{n-k-2} \cdot (-6)^{k+2} -(-3)^{n-k-2}\cdot 2^{k+2}}{2}.
\end{eqnarray*}
Here we used Lemma \ref{l-V1-p=2-HS-compute} to compute (i) and (ii). 
We obtain 
\begin{eqnarray*}
&&  \frac{5^{n-k-1} \cdot (-6)^{k+1} -(-3)^{n-k-1}\cdot 2^{k+1}}{2} + \frac{3}{2} \cdot \frac{5^{n-k-2} \cdot (-6)^{k+2} -(-3)^{n-k-2}\cdot 2^{k+2}}{2}\\
&=&  \frac{5^{n-k-1} \cdot (-6)^{k+1} }{2} + \frac{3}{2} \cdot \frac{5^{n-k-2} \cdot (-6)^{k+2}}{2} \\
&=& \frac{5  + \frac{3}{2} \cdot (-6)}{2} 
\cdot  5^{n-k-2}\cdot (-6)^{k+1}  
\\
&=& 12 \cdot 5^{n-2} \cdot (- \frac{6}{5})^k. 
\end{eqnarray*}
Therefore, 
\begin{eqnarray*}
 &&(4\tilde{H}^2+6\tilde{H}S +2S^2) \sum_{k=0}^{\infty} (k+1) \binom{n}{n-k-2} (\tilde{H}+2S)^{n-k-2}(-(2\tilde{H}+2S))^k\\
 &=& {\rm (i)} + \frac{3}{2} \times {\rm (ii)}\\
&=& \sum_{k=0}^{\infty} (k+1) \binom{n}{n-k-2} (12 \cdot 5^{n-2} \cdot (- \frac{6}{5})^k ) \\
&=& 12 \cdot 5^{n-2} \cdot \sum_{k=0}^{\infty} (k+1) \binom{n}{n-k-2}  (- \frac{6}{5})^k  \\
&\overset{(\alpha)}{=}& 12 \cdot 5^{n-2}  \cdot f(-\frac{6}{5})\\
&\overset{(\beta)}{=}& 12 \cdot 5^{n-2} \frac{n(-\frac{6}{5}) ((-\frac{6}{5})+1)^{n-1} -((-\frac{6}{5})+1)^n +1}{(-\frac{6}{5})^2}  \\
&=& 12 \cdot  \frac{n(-6) ((-6)+5)^{n-1} -((-6)+5)^n +5^n}{36}  \\
&=& \frac{1}{3}(  (-1)^n(6n-1) + 5^n), 
\end{eqnarray*}
($\alpha$) and ($\beta$) follow from Step  \ref{s3-Chern-V1-p=2}.  
This completes the proof of Step \ref{s4-Chern-V1-p=2}. 
\end{proof}
By Step \ref{s1-Chern-V1-p=2}, Step \ref{s2-Chern-V1-p=2}, and 
Step \ref{s4-Chern-V1-p=2}, we obtain 
\begin{eqnarray*}
c_n(\Omega_{\tilde Y}^1 \otimes M) &=&4n (-1)^n  + 6n (-1)^{n-1} + \frac{1}{3}(  (-1)^n(6n-1) + 5^n)\\
&=&  \frac{1}{3} ( 5^n-(-1)^n).
\end{eqnarray*}
This complete the proof of Theorem \ref{t-Chern-V1-p=2}. 
\end{proof}

\begin{cor}\label{c-Chern-V1-p=2}
We use Notation \ref{notation:deg=1}. 
Then $\phi\colon X\to Y$ is separable, i.e.,
Theorem \ref{thm:embedding(d=1)} (5) holds.
\end{cor}
\begin{proof}
We may assume $p=2$.
Since $\alpha \colon \tilde{X} \to X$ and $\beta\colon \tilde{Y} \to Y$ 
are birational, it suffices to show show that $\tilde{\psi} \colon \tilde{X} \to \tilde Y$ is separable. 

By Proposition \ref{p-tilde-L}, we have $K_{\tilde Y} \sim \tilde{\phi}^*(K_{\widetilde Y} + \tilde B)$ for $\tilde B \coloneqq 3\pi^*H_V + 2F$.
By \cite[Proposition 0.1.2]{CD89} and \cite[Lemma 3.2]{Ful98}, if $c_n(\Omega^1_{\tilde Y} \otimes \sO_{\tilde Y}(2\tilde{B})) \neq 0$, then $\tilde{\phi}$ is separable.
Set $M\coloneqq 2\tilde{B}=2(3\pi^*H_V + 2F) = 6\tilde{H}+4S_{-}$. 
Then, by Theorem \ref{t-Chern-V1-p=2}, we obtain $c_n(\Omega^1_{\tilde Y} \otimes \sO_{\tilde Y}(M))= \frac{1}{3}(5^n-(-1)^n)\neq 0$, and the assertion holds.
\end{proof}

%% file: section4.tex
\section{Singularities of ladders}\label{sec:normality of ladder}

In this section, we investigate singularities of general ladders of smooth del Pezzo varieties.

\begin{prop}\label{prop:degree geq 3}
    Let $(X,L)$ be an $n$-dimensional smooth del Pezzo variety. 
    Set $d\coloneqq L^{n}$. 
    Let
     \[
   (X,L)=(X_n,L_n) \supset (X_{n-1},L_{n-1}) \supset \cdots \supset (X_1,L_1)
     \] be a general ladder.
    Suppose that one of the following holds:
    \begin{enumerate}
        \item $d\geq 3$.
        \item $d=2$ and $p\geq 3$. 
    \end{enumerate}
    Then $X_i$ is smooth for every $i>0$.  
\end{prop}

\begin{proof}
By induction, it is enough to show that $X_{n-1}$ is smooth (Theorem \ref{thm:ladder}). 

    Suppose that (1) holds. In this case, $L$ is very ample by Proposition \ref{prop:description}(1), and thus a general member $X_{n-1}\in |L|$ is smooth.

    Suppose that (2) holds. 
    In this case, the morphism $\phi_{|L|} \colon X \to \P^n$ induced by the complete linear system $|L|$ is a finite double cover by Theorem \ref{thm:embedding(d=2)}(2). 
    Since $p\geq 3$, it follows from \cite[Corollary 4.3]{Spr98} that a general member $X_{n-1}$ of $|L|$ is smooth. 
\end{proof}



\begin{lem}\label{lem:degree 2}
    Let $(X,L)$ be an $n$-dimensional smooth del Pezzo variety 
    such that $L^{n}=2$.
    Let
     \[
   (X,L)=(X_n,L_n) \supset (X_{n-1},L_{n-1}) \supset \cdots \supset (X_1,L_1)
     \] be a general ladder.
    Then $X_1$ is smooth.
\end{lem}
\begin{proof}
We may assume that $p=2$ by Proposition \ref{prop:degree geq 3}.
We use the same notation as in Theorem \ref{thm:embedding(d=2)}.
Let $H_0, H_1,\ldots,H_{n-1} \in |L|$ be general members.
Since $L=\phi_{|L|}^{*}\sO_{\P^{n}}(1)$ and $\phi_{|L|} : X \to \P^n$ 
is a finite separable morphism of degree $2$, the intersection $ H_0 \cap H_1 \cap \cdots \cap H_{n-1}$ consists of two closed points $P$ and $Q$ with $P \neq Q$. 
Take the blowup 
\[
\sigma \colon Z\to X
\]
at $P$ and $Q$, 
which coincides with the resolution of indeterminacies of 
the 
rational map 
$X \dashrightarrow \P^{n-1}$ induced by 
the linear system generated by
 $H_0, H_1,\ldots, H_{n-1}$. 
 Hence we obtain the morphism 
\[
\pi\colon Z\to \P^{n-1}.
\]
Since $X_1$ is smooth at $P$ and $Q$, 
$X_1$ is (isomorphic to) a general fiber of $\pi : Z \to \P^{n-1}$. 
Let 
\[
\kappa \coloneqq k(t_0, ..., t_{n-2})
\]
be the function field of $\P^{n-1}$, 
where $k(t_0,\ldots, t_{n-2}) = K(\mathbb A^{n-1})$ and 
$t_i \coloneqq x_i/x_{n-1}$. 
For its algebraic closure $\overline \kappa$, 
let $Z_{\kappa}$ (resp.~$Z_{\overline{\kappa}}$) be the generic fiber (the geometric generic fiber) of $\pi$.

\setcounter{step}{0}

\begin{step}\label{s1-degree 2}
Assume that $Z_{\overline{\kappa}}$ is not smooth (i.e., $X_1$ is not smooth).  
Then each of $Z_{\overline{\kappa}}$ and $X_1$ 
is a cuspidal cubic curve. 

\end{step}

\begin{proof}[Proof of Step \ref{s1-degree 2}]
Since $X_1$ (i.e., a general fiber of $\pi$) is a curve of genus one, 
so is the geometric generic fiber $Z_{\overline{\kappa}}$. 
In particular, the geometric generic fiber $Z_{\overline{\kappa}}$ is integral. 
Since the normalization $Z_{\overline{\kappa}}^N \to Z_{\overline{\kappa}}$  
of $Z_{\overline{\kappa}}$ is a universal homeomorphism \cite[Lemma 2.2(2)]{Tan18b},  
also the normalisation $X_1^N \to X_1$ of a general fiber $X_1$ of $\pi$ is a universal homeomorphism. 
In particular, each of $Z_{\overline{\kappa}}$ and $X_1$  is a cuspical cubic curve. 
This completes the proof of Step \ref{s1-degree 2}. 
\end{proof}

By Theorem \ref{thm:embedding(d=2)}(1), we can write 
\[
X=\{f=y^2+f_2(x_0,\ldots,x_{n})y+f_4(x_0,\ldots,x_{n})=0 \subset \P(1,\ldots,1,2)_{[x_0:\cdots:x_{n}:y]} \}, 
\]
where  $f_2(x_0, ..., x_n)$ and $f_4(x_0, ..., x_n)$ are homogeneous polynomials of degree $2$ and $4$, respectively. 
Since $\phi$ is separable  (Theorem \ref{thm:embedding(d=2)}(3)), we have $f_2(x_0, ..., x_n) \neq 0$. 

\begin{step}\label{s2-degree 2}
If $f_2(x_0, ..., x_n)$ contains a term $x_{i_0}x_{j_0}$ with 
$1 \leq i_0 < j_0 \leq n$, then $X_1$ is smooth.
\end{step}

\begin{proof}[Proof of Step \ref{s2-degree 2}]
We have 
\[
f_2(x_0, ..., x_n) = \sum_{i=0}^{n} a_ix_i^2 + \sum_{0 \leq i< j \leq n} b_{ij} x_ix_j
\]
for some $a_i, b_{ij} \in k$. By assumption, $b_{i_0j_0} \neq 0$. 
We consider an evaluation as follows:
\[
x_0 = s_0 x_{n-1} +t_0x_n, \quad 
x_1 = s_1x_{n-1} +t_1 x_n, 
\ldots, 
x_{n-2} = s_{n-2}x_{n-1} +t_{n-2} x_n. 
\]
Then, after this evaluation, the coefficient of $x_{n-1}x_n$ in 
$f_2$ will become 
\[
\sum_{0 \leq i<j \leq n} b_{ij} (s_it_j+t_is_j). 
\]
Since $b_{i_0j_0} \neq 0$, this is nonzero when  $s_i \in k$ and $t_i \in k$ are general. 
Since the ladder in the statement is assumed to be general, 
$X_1$ is $F$-split by Fedder's criterion.
Since a cuspidal cubic is not $F$-split, $X_1$ have to be smooth 
(Step \ref{s1-degree 2}). 
This completes the proof of Step \ref{s2-degree 2}. 
\end{proof}

By Step 2, we may assume that 
$f_2(x_0, ..., x_n) =\sum_{i=0}^{n} a_ix_i^2$ with $a_i \in k$. 
Replacing 
$\sum_{i=0}^{n} \sqrt{a_i}x_i$ by $x_n$, 
the problem is reduced  to the case when 
\[
f_2(x_0, ..., x_n)=x_n^2.
\]

\begin{step}\label{s3-degree 2}
$X_1$ is smooth.
\end{step}

\begin{proof}[Proof of Step \ref{s3-degree 2}]
We first prove the following claim.

\begin{claim*}
There exist $s_0, ..., s_{n-1}, t_0, ..., t_{n-1} \in k$ such that 
the coefficient of $x_{n-1}^3 x_n$ in the polynomial 
\[
f_4(s_0 x_{n-1} +t_0x_n, s_1x_{n-1} +t_1 x_n, \ldots, 
s_{n-2}x_{n-1} +t_{n-2} x_n, x_{n-1}, x_n)  \in k[x_{n-1}, x_n]
\]
is nonzero. 
\end{claim*}

\begin{proof}[Proof of Claim]
An arbitrary monomial of $f_4$ can be written as 
$x_i^4$, $x_i^3x_j$, $x_i^2x_j^2$, $x_i^2x_jx_k$, or
 $x_ix_jx_kx_{\ell}$, where $i, j, k, \ell$ are distinct. 

The polynomial $f_4(s_0 x_{n-1} +t_0x_n, s_1x_{n-1} +t_1 x_n, \cdots, 
s_{n-2}x_{n-1} +t_{n-2} x_n, x_{n-1}, x_n)$ is obtained by applying the following  evaluation:
\[
x_0 = s_0 x_{n-1} +t_0x_n, \quad 
x_1 = s_1x_{n-1} +t_1 x_n, \cdots, 
x_{n-2} = s_{n-2}x_{n-1} +t_{n-2} x_n. 
\]
Then we have
\begin{align*}
    x_0^3x_1 &=(s_0 x_{n-1} +t_0x_n)^3(s_1 x_{n-1} +t_1x_n)\\
    &=(s_0^3 t_1 +s_0^2t_0s_1)x_{n-1}^3x_n +(\text{the other terms}).
\end{align*}
\begin{align*}
x_0^2x_1x_2 &=(s_0 x_{n-1} +t_0x_n)^2(s_1 x_{n-1} +t_1x_n)(s_2 x_{n-1} +t_2x_n)\\ &=
(s_0^2t_1s_2+s_0^2s_1t_2)x_{n-1}^3x_n +(\text{the other terms}). 
\end{align*}
\begin{align*}
    x_0x_1x_2x_3 &=(s_0 x_{n-1} +t_0x_n)(s_1 x_{n-1} +t_1x_n)(s_2 x_{n-1} +t_2x_n)(s_3 x_{n-1} +t_3x_n)\\ 
&=(t_0s_1s_2s_3 + t_1s_0s_2s_3 + t_2s_0s_1s_3 +t_3s_0s_1s_2)x_{n-1}^3x_n +(\text{the other terms}).
\end{align*}
Note that the coefficient of $x_{n-1}^3 x_n$ in $x_0^4$ and $x_0^2 x_1^2$ is zero since we are working in characteristic two. 
We can write 
\begin{align*}
    &f_4(s_0x_{n-1}+t_0x_n, \ldots, s_{n-2}x_{n-1}+t_{n-2}x_{n-2}, x_{n-1}, x_n)\\
    =&\varphi(s_0, \ldots, s_{n-2}, t_0, \ldots, t_{n-2})x_{n-1}^3 x_n + (\text{the other terms}). 
\end{align*}
We treat (I) and (II) separately. 
\begin{enumerate}
\item[(I)]
$f_4(x_0, \ldots, x_n)$ does not have any of the monomials $x_i^3x_j$, $x_i^2x_jx_k$, and $x_ix_jx_kx_{\ell}$. 
\item[(II)] $f_4(x_0,\ldots, x_n)$ has one of the monomials $x_i^3x_j$, $x_i^2x_jx_k$, and $x_ix_jx_kx_{\ell}$. 
\end{enumerate}

(I) In this case, we derive a contradiction by showing that $X$ is not smooth. 
We can write 
\[
f_4(x_0, \ldots, x_n) = \sum_{i} \alpha_i x_i^4 + \sum_{i, j} \beta_{ij} x_i^2x_j^2. 
\]
Then 
\[
X \cap D_+(x_0) = \{f =  y^2 + x_n^2 y + f_4(1, x_1, \ldots, x_n) =0\} \in \A^{n+1}_{(x_1, \ldots, x_n, y)}. 
\]
Let us find a singular point on $X \cap D_+(x_0)$ by Jacobian criterion. 
Clearly, we have 
\[
\partial_{x_1}f = \cdots = \partial_{x_n}f =0, 
\]
where each $\partial_{x_i}$ denotes the partial differential. 
Then it is enough to find a solution of 
\[
f = \partial_y f=0. 
\]
For $b \in k$ satisfying $ b^2 + f_4(1, 0,\ldots, 0)=0$, 
we see that $(x_1, \ldots, x_n, y) = (0, \ldots, 0, b)$ is a solution of this equation. 
This is absurd. 

\medskip

(II) 
In order to prove Claim, it is enough to show that 
the coefficient 
$\varphi(s_0, ..., s_{n-2}, t_0, ..., t_{n-2})$ 
of $x_{n-1}^3x_n$ is nonzero when $s_i \in k$ and $t_i \in k$ are general. 
To this end, we consider $s_0, ... s_{n-2}, t_0, ..., t_{n-2}$ as variables, and it suffices check that the coefficient $\varphi(s_0, ..., s_{n-2}, t_0, ..., t_{n-2})$  of $x_{n-1}x_n$ is nonzero. 
In this situation, we have 
\[
\varphi(s_0, ..., s_{n-2}, t_0, ..., t_{n-2}) \in k[s_0, ..., s_{n-2}, t_0, ..., t_{n-2}]. 
\]
If $f_4(x_0, ..., x_n)$ has $x_0^3x_1$, 
then the monomial $s_0^3 t_1$ appears in 
$\varphi(s_0, ..., s_{n-2}, t_0, ..., t_{n-2})$, because 
$s_0^3t_1$ does not appear from the other terms. 
By symmetry, we may assume that any monomial of the form $x_i^3x_j$ does not appear in $f_4(x_0, ..., x_n)$ . 
If $f_4(x_0, ..., x_n)$ has $x_0^2x_1x_2$, 
then the monomial $s_0^2t_1s_2$ is alive in $\varphi(s_0, ..., s_{n-2}, t_0, ..., t_{n-2})$. Hence we may assume that 
$f_4(x_0, ..., x_n)$ has no monomial of the form $x_i^2x_jx_k$. 
Similarly, we may exclude the case when 
$f_4(x_0, ..., x_n)$ has a  monomial of the form $x_ix_jx_kx_{\ell}$. 
This  finishes the proof of Claim.
\qedhere


\end{proof}

Take general $s_0, ..., s_{n-1}, t_0, ..., t_{n-1} \in k$ as in Claim. 
Set 
\[
g_4(x_{n-1}, x_n) := f_4(s_0 x_{n-1} +t_0x_n, s_1x_{n-1} +t_1 x_n, \cdots, 
s_{n-2}x_{n-1} +t_{n-2} x_n, x_{n-1}, x_n). 
\]
By Claim, $g_4(x_{n-1}, x_n)\in k[x_{n-1}, x_n]$ is a homogeneous polynomial of degree $4$ which has the term $x_{n-1}^3x_n$. 
We may assume that $X_1$ is given by 
\[
X_1 = \{ [x_{n-1}:x_n : y] \,|\, y^2 + x_n^2 y +g_4(x_{n-1}, x_n) =0\} 
\subset \P(1, 1, 2)_{[x_{n-1}:x_n : y]}. 
\]
Since $X$ does not pass through $[0:\cdots:0:1] \in \P(1, ..., 1, 1, 2)_{[x_0: \cdots : x_{n-1}:x_n:y]}$, it follows that $X_1$ does not pass through 
\[
[0:0:1] \in \P(1, 1, 2)_{[x_{n-1}:x_n:y]} 
\]
\[
= \P(1, ..., 1, 1, 2)_{[x_0: \cdots : x_{n-1}:x_n:y]} \cap 
\left(\bigcap_{0 \leq i\leq n-2} \{x_i = s_i x_{n-1} +t_ix_n\}\right). 
\]
Thus it suffices to show that $X_1$ is smooth on 
each of $D_{+}(x_n)$ and $D_{+}(x_{n-1})$.

On $D_{+}(x_n)$, we have
\[
X_{1}|_{D_{+}(x_n)}=\{h(x_{n-1},y):=y^2+y+ g_4(x_{n-1},1)=0\} \in \A^2_{(x_{n-1},y)}. 
\]
Since $\partial_{y}h=1$, the Jacobian criterion shows that $X_{1}|_{D_{+}(x_n)}$ is smooth.

On $D_{+}(x_{n-1})$, we have
\[
X_{1}|_{D_{+}(x_{n-1})}=\{\wt{h}(x_{n},y):=y^2+x_{n}^2y+a_4 x_{n}^4+a_3 x_{n}^3 + a_2 x_{n}^2 + a_{1} x_{n} + a_0=0\} \in \A^2_{(x_n,y)}
\]
with $a_i \in k$. 
Since $g_4(x_{n-1}, x_n)$ has the term $x_{n-1}^3x_n$, 
we get $a_1\neq 0$.
Since $\partial_{x_n}h=a_3x_{n}^2+a_{1}$ and $\partial_{y}h=x_{n}$, the Jacobian criterion shows that $X_{1}|_{D_{+}(x_{n-1})}$ is smooth.
This completes the proof of Step \ref{s3-degree 2}. 
\end{proof}
Step \ref{s3-degree 2} compltes the proof of Lemma \ref{lem:degree 2}. 
\end{proof}

\begin{lem}\label{l cusp triple}
Assume $p=3$. 
Let 
\[
C := \{ (y, z) \in \A^2\,|\, z^2 + y^3 + a y^2 + by + c =0\} 
\]
be an affine curve, where $a, b, c \in k$. 
Assume that $C$ contains a singularity which is analytically isomorphic to 
the cuspidal singularity $k[X, Y]_{(X, Y)}/(Y^2 + X^3)$. 
Then $a = b=0$. 
\end{lem}

\begin{proof}
By taking the closure $\overline{C} := \{ [x:y:z] \in \P^2 \,|\, xz^2 + y^3 + a xy^2 + bx^2y + cx^3 =0\}$, 
we see that $C$ has a unique singularity $Q$, 
as otherwise we would get a contradiction $3 = \overline{C} \cdot L \geq 2 + 2$ 
for the line $L$ passing through two of the singular points of $\overline C$.  
By the Jacobian criterion, the  singular point $Q$ of $C$ is a solution of 
\[
z =  y^3 + a y^2 + by + c  = 
\frac{d}{dy}(y^3 + a y^2 + by + c)=0. 
\]
Then we have  $Q = (y_1, 0) \in k^2$, and 
$y=y_1$ is a multiple solution of $y^3 + ay^2 + by +c=0$. 
Hence we can write 
\[
y^3 + ay^2 + by + c = (y-y_1)^2 (y-y_2)
\]
for some $y_2 \in k$. 
If $y_1 \neq y_2$, then $Q=(y_1, 0)$ would be a nodal singularity, because 
$z^2 = (y-y_1)^2(y-y_2)$ will become $z^2 = y^2(y+1) = y^2 + y^3$ after a suitable coordinate change. 
Therefore, we get $ y_1 = y_2$, which implies 
\[
y^3 + ay^2 + by + c = (y-y_1)^3 = y^3 - y_1^3. 
\]
Hence $a =b=0$.
\end{proof}

\begin{lem}\label{lem:degree 1}
    Let $(X,L)$ be an $n$-dimensional smooth del Pezzo variety 
    such that $L^{n}=1$.
        Let
     \[
   (X,L)=(X_n,L_n) \supset (X_{n-1},L_{n-1}) \supset \cdots \supset (X_1,L_1)
     \] be a general ladder of $(X, L)$. 
    Then 
    $X_1$ is smooth. 
\end{lem}

\begin{proof}
We use the same notation as in  Theorem \ref{thm:embedding(d=1)}.
As in Lemma \ref{lem:degree 2}, we can take effective Cartier divisors 
$H_0, H_1,\ldots,H_{n-1} \in |L|$ such that $\bigcap_{i=1}^{n-1}H_{i}=X_{1}$ 
and $P  := \Bs |L| = H_0 \cap H_1 \cap \cdots \cap H_{n-1}$. 
Take the blowup 
\[
\sigma : Z\to X
\]
at $P=\Bs |L|= H_0 \cap H_1 \cap \cdots \cap H_{n-1}$. 
Then  we obtain the induced morphism 
\[
\pi\colon Z\to \P^{n-1}, 
\]
because  $\sigma : Z \to X$ coincides with the resolution of indeterminacies of 
the rational map $X \dashrightarrow \P^{n-1}$ induced by 
the complete linear system $|L|$, which is generated by
 $H_0, H_1, ..., H_{n-1}$. 
 We may assume that $H_i = \{ x_i =0\}$ for every $0 \leq i \leq n-1$. 
Since $X_1$ is smooth at $P$ (Lemma \ref{l deg1 regularity}), 
$X_1$ is isomorphic to its proper transform $X_1^Z$ on $Z$, 
which is a general fiber of $\pi$. 
In what follows, we shall identify $X_1$ and $X_1^Z$. 
Let 
\[
\kappa := k(t_1, ..., t_{n-1})
\]
be the function field of $\P^{n-1}$, 
where $k(t_1, ..., t_{n-1}) = K(\mathbb A^{n-1})$ and $t_i := x_i/x_0$. 
For its algebraic closure $\overline \kappa$, 
let $Z_{\kappa}$ (resp. $Z_{\overline{\kappa}}$) be the generic fiber (the geometric generic fiber) of $\pi : Z \to \P^{n-1}$. 
We have 
\[
X = \{ [x_0 : \cdots :x_{n-1} : y :z ] \in \P(1, ..., 1, 2, 3) \,|\, f(x_0, ..., x_{n-1}, y, z)=0\}
\]
for a homogeneous polynomial $f(x_0, ..., x_{n-1}, y, z)$ of degree $6$. 
For a general closed point $[a_0: \cdots : a_{n-1}] \in \P^{n-1}$ satisfying $a_0 =1$, it holds that 
\begin{eqnarray*}
&&\pi^{-1}([a_0: \cdots : a_{n-1}]) \\
&\xrightarrow{\sigma, \simeq}& 
X \cap \left(\bigcap_{0 \leq i<j\leq n-1} \{ [x_0 : \cdots :x_{n-1} :y:z] \in X\,|\, 
a_jx_i = a_i x_j\}\right)\\
&=& X \cap \left(\bigcap_{j =1}^{n-1} \{ [x_0 : \cdots :x_{n-1} :y:z] \in X\,|\, 
a_jx_0 = a_0 x_j\}\right)\\
&=& \{ [x_0:y:z] \in \P(1, 2, 3)_{[x_0:y:z]} \,|\, 
f(x_0, a_1x_0, ...,  , a_{n-1}x_0, y, z)\}.
\end{eqnarray*}
Similarly, the generic fiber $Z_{\kappa}$ of $\pi$ satisfies 
\begin{eqnarray*}
Z_{\kappa} &\simeq& 
\{ [x_0:y:z] \in \P(1, 2, 3)_{[x_0:y:z]} \times_k \kappa \,|\, 
f(x_0, t_1x_0, ...,  , t_{n-1}x_0, y, z) =0\}\\
&=& \Proj\, \frac{\kappa[x_0, y, z]}{f(x_0, t_1x_0, ...,  , t_{n-1}x_0, y, z)}.  
\end{eqnarray*}

\setcounter{step}{0}

\begin{step}\label{s1-degree 1}
Assume that $Z_{\overline{\kappa}}$ is not smooth (i.e., $X_1$ is not smooth).  
Then each of $Z_{\overline{\kappa}}$ and $X_1$ 
is a cuspidal cubic curve. 

\end{step}

\begin{proof}[Proof of Step \ref{s1-degree 1}]
Since $X_1$ (i.e., a general fiber of $\pi$) is a curve of genus one, 
so is the geometric generic fiber $Z_{\overline{\kappa}}$. 
In particular, the geometric generic fiber $Z_{\overline{\kappa}}$ is integral. 
Since the normalization $Z_{\overline{\kappa}}^N \to Z_{\overline{\kappa}}$  
of $Z_{\overline{\kappa}}$ is a universal homeomorphism \cite[Lemma 2.2(2)]{Tan18b},  
also the normalisation $X_1^N \to X_1$ of a general fiber $X_1$ of $\pi$ is a universal homeomorphism. 
In particular, each of $Z_{\overline{\kappa}}$ and $X_1$  is a cuspical cubic curve. 
This completes the proof of Step \ref{s1-degree 1}. 
\end{proof}

\begin{step}\label{s2-degree 1}
If $p \geq 5$, then $X_1$ is smooth. 
\end{step}

\begin{proof}[Proof of Step \ref{s2-degree 1}]
Since $Z_{\kappa}$ is a regular projective curve with $h^1(Z_{\kappa}, \MO_{\kappa})=1$,    
it follows from $p \geq 5$ and \cite[Theorem 1.4]{PW22} that $Z_{\kappa}$ is smooth over $\kappa$. 
Hence a general fiber $X_1$ of $\pi : Z \to \P^{n-1}$ is smooth. 
This completes the proof of Step \ref{s2-degree 1}. 
\end{proof}

\begin{step}\label{s3-degree 1}
If $p =3$, then $X_1$ is smooth. 
\end{step}

\begin{proof}[Proof of Step \ref{s3-degree 1}]
Assume $p=3$.
We can write 
\[
X=\{z^2+(f_3
+yf_1 )z+
y^3+y^2f_2
+yf_4
+f_6
=0 \}\subset \P(1,\ldots,1,2,3)_{[x_0:\cdots:x_{n-1}:y:z]}, 
\]
where each $f_i := f_i(x_0,\ldots,x_{n-1}) \in k[x_0, ..., x_{n-1}]$ 
is a homogeneous polynomial of degree $i$. 
By applying  $z\mapsto z-(1/2)(f_3 + yf_1)$, 
we may assume that $f_1 = f_3=0$. 
Hence 
\[
X=\{z^2+
y^3+y^2f_2
+yf_4
+f_6
=0 \}\subset \P(1,\ldots,1,2,3)_{[x_0:\cdots:x_{n-1}:y:z]}. 
\]

\begin{claim*}
$f_2 = f_4=0$. 
\end{claim*}

\begin{proof}[Proof of Claim]
We only prove $f_2 =0$, as both proofs are identical. 
Suppose $f_2 = f_2(x_0, ..., x_{n-1})\neq 0$. 
It suffices to derive a contradiction. 
Pick a general closed point $[a_0: \cdots : a_{n-1}] \in\P^{n-1}_k$. 
By Step \ref{s1-degree 1}, we may assume that 
\begin{itemize}
\item[(i)] 
$a_0\neq 0, a_1 \neq 0, ..., a_{n-1} \neq 0$, 
\item[(ii)] 
$f_2(a_0, \cdots, a_{n-1}) \neq 0$, and 
\item[(iii)] 
the singular point of the cuspidal cubic curve 
\[
F := \pi^{-1}([a_0 : \cdots : a_{n-1}]) =  \{ z^2 + y^3 +b_2 x_0^2y^2 + b_4 x_0^4y 
+ b_6 x_0^6 =0\} \subset \P(1, 2, 3)_{[x_0:y:z]} 
\]
is contained in $D_+(x_0)$, 
where $b_i := f_i(a_0, ..., a_{n-1}) \in k$ for each $i \in \{2, 4, 6\}$. 
\end{itemize}
We now verify how to assure (iii). 
Recall that 
\[
X \cap \{x_0 = \cdots = x_{n-1}=0\} =H_0 \cap H_1 \cap \cdots \cap H_{n-1} = P. 
\]
Then a general member $X_1$ is smooth at $P$ (Lemma \ref{l deg1 regularity}(3)). 
Therefore, also $\sigma(F) (\simeq F)$ is smooth at $P$. 
Hence the singular point of $\sigma(F)$ is contained in $\bigcup_{i=0}^{n-1} D_+(x_i)$. 
After permuting the indices, the problem is reduced to the case when (i)-(iii) hold.

 By (i), we  may assume that $a_0 =1$. 
For each $i \in \{2, 4, 6\}$,  we have 
\[
f_i(x_0, a_1x_0, ..., a_{n-1}x_0)= f_i(a_0, a_1, ..., a_{n-1}) x_0^i = b_i x_0^i. 
\]
It follows from (iii) that the affine plane cubic curve 
\[
F \cap D_+(x_0) = \{ (y, z) \in \A^2 \,|\, z^2 +y^3 + b_2 y^2 + b_4y + b_6=0\}
\]
contains a cuspidal singularity. 
We then get $b_2 = b_4 =0$ (Lemm \ref{l cusp triple}), which is  absurd, 
because (ii) and (iii) imply  $b_2= f_2(a_0, a_1, ..., a_{n-1})\neq 0$. 
This completes the proof of Claim. 
\end{proof}

\medskip

By Claim, the problem has been reduced to the case when 
\[
X=\{z^2+
y^3+f_6
=0 \}\subset \P(1,\ldots,1,2,3)_{[x_0:\cdots:x_{n-1}:y:z]}. 
\]
We shall prove  that 
\begin{enumerate}
\item[$(\star)$] 
$Y\coloneqq \{y^3+f_6(x_0,\ldots,x_{n-1})=0 \}\subset \P(1,\ldots,1,2)_{[x_0:\cdots:x_{n-1}:y]}$ is smooth. 
\end{enumerate}
We now finish the proof of Step \ref{s3-degree 1} by assuming ($\star$). 
We have a finite inseparable morphism of degree $3$ 
\[
Y\to \P^{n-1};\qquad [x_0 :\cdots : x_{n-1}:y] \mapsto [x_0 :\cdots : x_{n-1}].
\]
This contradicts ($\star$), because a smooth inseparable cover of $\P^n$ of degree $p$ 
does not exist in characteristic $p \geq 3$ \cite[Proposition 2.5]{Eke87}.


It is enough to show $(\star)$. 
Suppose that $Y$ has a singular point $Q := [a_0: \cdots : a_{n-1}:b] \in Y \subset \P(1, ..., 1, 2)$ with $a_0, ..., a_{n-1}, b \in k$. 
It suffices to derive a contradiction. 
If $a_0= \cdots =a_{n-1}=0$, then the defining equation $y^3 + f_6(x_0, ..., x_{n-1})=0$ would imply $b=0$, which is absurd. 
Hence one of $a_0, ..., a_{n-1}$ is nonzero. 
After permuting the indices, 
we may assume that $a_0 =1$. 
In this case,  $Q \in D_+(x_0) \subset \P(1, ..., 1, 2)$ is 
a singular point of the affine open subset $Y'$ of $Y$ given as follows: 
\[
Y' = \{ (\wt{x}_1,\ldots, \wt{x}_{n-1}, \wt{y})\,|\, g(\wt{x}_1,\ldots, \wt{x}_{n-1}, \wt{y}) =0\} \subset \Spec k \left[ \wt{x}_1, \wt{x}_2, ..., \wt{x}_{n-1}, \wt{y}\right] = D_+(x_0)
\]
\[
g(\wt{x}_1, ..., \wt{x}_{n-1}, \wt{y}) := \wt{y}^3 +f_6(1, \wt{x}_1, ..., \wt{x}_{n-1}), 
\]
where $\wt{x}_i := x_i/x_0$ and $\wt{y}:= y/x_0^2$. 
By the Jacobian criterion, 
$(a_1, ..., a_{n-1}, b) \in k^n$ (which is a singular point of $Y'$) 
is a solution of  
\[
g= \partial_{\wt{x}_1} g = \cdots =\partial_{\wt{x}_{n-1}}g = \partial_{\wt{y}} g =0, 
\]
where each $\partial_{\bullet}$ denotes the partial differential. 
It is enough to show that $P := [a_0 : a_1 : \cdots : a_{n-1}:b:0] \in \P(1, ..., 1, 2, 3)$ is a singular point of $X$. 
Take the affine open subset $D_+(x_0) := \{ x_0 \neq 0\}$ of $\P(1, ..., 1, 2, 3)$: 
\[
X' \coloneqq \{ h(\wt{x}_1,\ldots, \wt{x}_{n-1}, \wt{y}, \wt{z}) =0\} \subset \Spec k[\wt{x}_1, \ldots, \wt{x}_{n-1}, \wt{y}, \wt{z}] = D_+(x_0)
\]
\[
h(\wt{x}_1, \ldots, \wt{x}_{n-1}, \wt{y}, \wt{z}) \coloneqq \wt{z}^2 + \wt{y}^3 + f_6(1, \wt{x}_1,\ldots, \wt{x}_{n-1}), 
\]
where $\wt{x}_i \coloneqq x_i/x_0, \wt{y}\coloneqq y/x_0^2, \wt{z}\coloneqq z/x_0^3$. 
It suffices to prove that $(\wt{x}_1, \ldots, \wt{x}_{n-1}, \wt{y}, \wt{z}) = (a_1,\ldots, a_{n-1}, b, 0)$ is a solution of 
\[
h = \partial_{\wt{x}_1} h = \cdots = \partial_{\wt{x}_{n-1}} h = \partial_{\wt{y}} h = \partial_{\wt{z}} h =0. 
\]
This follows from  
\[
h(\wt{x}_1,\ldots, \wt{x}_{n-1}, \wt{y}, 0)=
g(\wt{x}_1,\ldots, \wt{x}_{n-1}, \wt{y}), 
\]
\[
\partial_{\wt{x}_1} g =\partial_{\wt{x}_1} h,\,\,\ldots,\,\, \partial_{\wt{x}_{n-1}} g =\partial_{\wt{x}_{n-1}} h,\,\,
\partial_{\wt{y}} g =\partial_{\wt y} h,\,\, \partial_{\wt{z}}h = 2\wt{z}. 
\]
Thus $(\star)$ holds. 
This completes the proof of Step \ref{s3-degree 1}. 
\end{proof}

\begin{step}\label{s4-degree 1}
If $p =2$, then $X_1$ is smooth. 
\end{step}

\begin{proof}[Proof of Step \ref{s4-degree 1}]
Assume  $p=2$ and $X_1$ is not smooth. 
Let us derive a contradiction. 
By Step \ref{s1-degree 1}, $X_1$ is a cuspidal cubic curve. 
We can write 
\[
\{z^2+(f_3+yf_1)z+y^3 + f_2y^2 + yf_4 + f_6=0 \}\subset \P(1,\ldots,1,2,3)_{[x_0:\cdots:x_{n-1}:y:z]}.
\]
where each $f_i := f_i(x_0, ..., x_{n-1}) \in k[x_0, ..., x_{n-1}]$ is a homogeneous polynomial of degree $i$. 
Since the double cover 
\[
\phi : X \to \P(1, ..., 1, 2); \quad [x_0: \cdots :x_{n-1}:y:z] \mapsto 
[x_0: \cdots :x_{n-1}:y]
\]
is separable by Theorem \ref{thm:embedding(d=1)}, we get 
\[
f_3 + yf_1  \neq 0. 
\]
If $f_1 \neq 0$, then $X_1$ is $F$-split by Fedder's criterion. 
However, this contradicts the fact that $X_1$ is the cuspidal cubic curve, 
which is not $F$-split. 
Hence it holds that  
\[
f_1 =0\qquad \text{and} \qquad f_3 \neq 0. 
\]

Recall that 
we fix a general ladder as in the statement,  
\begin{itemize}
\item $X_1 = X \cap  \{ x_1 = \cdots =x_n=0\}$, and 
\item $a:=f_3(1, 0, ..., 0) \neq 0$. 
\end{itemize}
It holds that $X_1$ is smooth at $P=\Bs |H|=\{x_0=\cdots=x_{n-1}\}$ 
(Lemma \ref{l deg1 regularity}(3)), 
where the latter equality is guaranteed by the proof of Theorem \ref{thm:embedding(d=1)}(3). 
It is enough to show that $X_1|_{D_+(x_i)}$ is smooth for every $0 \leq i \leq n$. 
Since $X_1$ is contained in $\{x_1 = \cdots = x_n=0\}$, 
it suffices to prove that $X_1$ is smooth on $D_{+}(x_0)=\Spec\,k[x_1/x_0,\ldots,x_{n-1}/x_0,y/x_0^2,z/x_0^3]\cong \mathbb{A}^{n+1}$.
Substituting $x_0=1, x_1=0, ..., x_{n-1}=0$,  it suffices to show that 
the following affine plane curve 
\[
X_1|_{D_{+}(x_0)}=\{f=z^2+az+h(y)=0 \}\subset \mathbb{A}^{2}_{(y, z)}
\]
is smooth, 
where
$h(y) \coloneqq f_6(1, 0, ..., 0, y) \in k[y]$. 
Since $\partial f/\partial z=a\neq 0$, the Jacobican criterion implies
that $X_1$ is smooth.
This completes the proof of Step \ref{s4-degree 1}. 
\end{proof}
Step \ref{s2-degree 1}, Step \ref{s3-degree 1}, and Step \ref{s4-degree 1}  complete the proof of Lemma \ref{lem:degree 1}. 
\end{proof}

\begin{thm}\label{thm:normality of ladder}
    Let $(X,L)$ be an $n$-dimensional smooth del Pezzo variety 
    and
    let
    \[
   (X,L)=(X_n,L_n) \supset (X_{n-1},L_{n-1}) \supset \cdots \supset (X_1,L_1)
    \] be a general ladder of $(X, L)$. 
    Then $X_i$ is normal for all $i\geq 1$.
\end{thm}
\begin{proof}
We prove that $X_i$ is normal by induction on $i$. 
    By Proposition \ref{prop:degree geq 3}, and Lemmas \ref{lem:degree 2} and \ref{lem:degree 1}, it follows that $X_1$ is smooth.
Fix $i >1$ and assume that $X_{i-1}$ is normal. 
It is enough to show that $X_i$ is normal. 
        Suppose by contradiction that $X_{i}$ is non-normal.
        Since $X_i$ is Cohen-Macaulay, it follows that $X_i$ is not $(R_1)$, i.e., the singular locus $Z_i$ of $X_i$ is of codimension one.
        Then $X_{i-1}\cap Z_i\subset X_{i-1}$ has codimension at least one since $X_{i-1}$ is an ample effective Cartier divisor on $X_i$.
        This contradicts the fact that $X_{i-1}$ is $(R_1)$.
\end{proof}

\begin{thm}\label{thm:X_2 is a can. dP}
    In the notation of Theorem \ref{thm:normality of ladder}, $X_2$ is a canonical del Pezzo surface.
\end{thm}

\begin{proof}
We see that $X_1 = H_1 \cap \cdots \cap H_{n-1}$ for some $H_1, ..., H_{n-1} \in |L|$. 
Let $\sigma \colon Z \to X$ be the blowup along $X_1$, 
which coincides with the resolution of the indeterminacies 
of the rational map $X \dashrightarrow \P^{n-2}$ 
of the linear system generated by $H_1, ..., H_{n-1}$. 
We then get the induced morphism $\pi : Z \to \P^{n-2}$. 
Since $X_2$ is smooth along the blowup center $X_1$, 
we get $X_2^Z \xrightarrow{\sigma, \simeq} X_2$ 
for the proper transform $X_2^Z$ of $X_2$ on $Z$. 
Since $X_2^Z$ is a (general) fiber of $\pi$, 
the generic fiber $Z_{K} := Z \times_{\P^{n-2}} \Spec K$ of $\pi$ is a regular del Pezzo surface, 
where $K$ denotes the function field of $\P^{n-2}$. 
Since $X_2$ is normal (Theorem \ref{thm:normality of ladder}), 
$Z_{K}$ is geometrically normal. 
Therefore, the geometric generic fiber $Z_{\overline K}$ is a del Pezzo surface which has at worst canonical singularities (\cite[Theorem 3.3]{BT22}). 
Then a general fiber $X_2$ of $\pi$ also has at worst canonical singularities. 
\end{proof}

%% file: section5.tex
\section{Weak quasi-$F$-splitting}\label{sec:quasi-F-split}

The purpose of this section is to introduce weak quasi-$F$-splitting  (Subsections \ref{ss-Witt-div}
\ref{ss-Witt-FV}, \ref{s-weak-QFS}) 
and establish a log Calabi-Yau inversion of adjunction 
for weak quasi-$F$-splitting (Subsection \ref{ss-logCY-IOA}).

\subsection{Witt divisorial sheaves}\label{ss-Witt-div}
We shall use the following variant of Witt divisorial sheaves introduced in \cite[Definition 3.1]{Tan22}. 
Comparison to the original definition, the following one involves the ideal term $I_E$. 

\begin{dfn}\label{d-div-def}
Let $X$ be an integral normal noetherian $\F_p$-scheme. 
Let $D$ be an $\R$-divisor on $X$ and let $E$ be an effective $\R$-divisor on $X$. 
Then we define the subpresheaf $WI_{X, E}(D)$ (or $WI_E(D)$) of the constant sheaf $W(K(X))$ on $X$ 
by 
{\small 
$$\Gamma(U, WI_E(D)):=\{(\varphi_0, \varphi_1, \cdots) \in W(K(X))\,|\,
\left({\rm div}(\varphi_n)+p^nD -E\right)|_U \geq 0\text{ for all }n\geq 0\}$$
}
for any open subset $U$ of $X$, 
where ${\rm div}(\varphi_n)$ denotes the principal divisor associated to $\varphi_n$. 
By definition, $WI_E(D)$ is a subsheaf of $W(K(X))$ (cf.~Remark \ref{r-div-def}). 
We call $WI_E(D)$ the {\em Witt divisorial sheaf associated to} $(D, E)$. 
We define the subsheaf $W_nI_E(D)$ of the constant sheaf $W_n(K(X))$ 
in the same way. 
\end{dfn}

\begin{rem}\label{r-div-def}
We use notation of Definition \ref{d-div-def}. If we identify $W(K(X))$ with the infinite direct product $\prod_{n=0}^{\infty}K(X)$ as sets, then it follows by definition that as subsets of $W(K(X))=\prod_{n=0}^{\infty}K(X)$, we obtain an equation 
\[
\Gamma(U, WI_E(D))=\prod_{n=0}^{\infty} \Gamma(U, \MO_X(p^nD-E)) = \prod_{n=0}^{\infty} \Gamma(U, \MO_X(\llcorner p^nD-E\lrcorner)).
\]
In particular, we have 
\[
W_1I_{\Delta}(p^e \Delta) =\MO_X( \llcorner (p^e-1)\Delta \lrcorner), 
\]
which is what appears in the definition of $F$-purity. 
\end{rem}

\begin{lem}\label{l-div-sub}
Let $X$ be an integral normal noetherian $\F_p$-scheme. 
Let $D$ be an $\R$-divisor  and let $E$ be an effective $\R$-divisor on $X$. 
Then the following hold. 
\begin{enumerate}
\item 
$WI_E(D)$ is a sheaf of $W\MO_X$-submodules of $W(K(X))$, 
i.e., for any open subset $U$ of $X$, 
$\Gamma(U,WI_E(D))$ is a $W\MO_X(U)$-submodule of $W(K(X))$. 
\item 
For any positive integer $n$, 
$W_nI_E(D)$ is a sheaf of $W_n\MO_X$-submodules of $W_n(K(X))$, 
i.e., for any open subset $U$ of $X$, 
$\Gamma(U, W_nI_E(D))$ is a $W_n\MO_X(U)$-submodule of $W_n(K(X))$. 
Furthermore, $W_nI_E(D)$ is a quasi-coherent $W_n\MO_X$-module. 
\end{enumerate}
\end{lem}

\begin{proof}
The same argument as in \cite[Lemma 3.5]{Tan22} works. 
\end{proof}

\begin{rem}\label{r-div-codim2}
We use the same notation as in Definition \ref{d-div-def}. 
Take an open subset $X'$ of $X$ such that 
$X'$ contains all the points of codimension one in $X$. 
Let $j:X' \to X$ be the induced open immersion. 
Then it holds that both 
$WI_E(D) \to j_*(WI_E(D)|_{X'})$ and 
$W_nI_E(D) \to j_*(W_nI_E(D)|_{X'})$ 
are isomorphisms for any positive integer $n>0$ (cf.~Remark \ref{r-div-def}). 
\end{rem}

\subsection{Frobenius and Verschiebung}\label{ss-Witt-FV}

\begin{dfn}
Let $X$ be an integral normal noetherian $\F_p$-scheme and 
let $D$ be an $\R$-divisor on $X$ and let $E$ be an effective $\R$-divisor. 
\begin{enumerate}
\item $F^e:W_n(K(X)) \to F_*^eW_n(K(X))$ induces 
\[
F^e:W_nI_E(D) \to F_*^e W_nI_{p^eE}(p^eD) (\subset F_*^e W_nI_{E}(p^eD)).
\]
\item $V^e:F_*^eW_n(K(X)) \to W_{n+e}(K(X))$ induces 
\[
V^e:F_*^e (W_nI_E(D)) \to F_*^e W_nI_{E}((1/p^e)D). 
\]
\item $R^m:W_{n+m}(K(X)) \to W_n(K(X))$ induces 
\[
R^m: W_{n+m}I_E(D) \to W_nI_{E}(D). 
\]
\end{enumerate}
\end{dfn}

\begin{prop}\label{p-div-induction}
Let $X$ be an integral normal noetherian $\F_p$-scheme and 
let $D$ be an $\R$-divisor on $X$ and let $E$ be an effective $\R$-divisor.  
Then there are exact sequences of $W\MO_X$-module homomorphisms: 
\begin{equation}\label{e-div-induction1}
0 \to (F^n_X)_*(WI_E(p^nD) ) \xrightarrow{V^n} WI_E(D) \to W_nI_E(D)\to 0
\end{equation}
{\small 
\begin{equation}\label{e-div-induction2}
0 \to (F^n_X)_*(W_mI_E(p^nD) ) \xrightarrow{V^n} W_{n+m}I_E(D) \to W_nI_E(D)\to 0
\end{equation}
}
for any positive integers $n$ and $m$. 
\end{prop}

\begin{proof}
Let us prove that 
the first sequence (\ref{e-div-induction1}) is exact. 
This is a subcomplex of the exact sequence 
$$0 \to (F_X^n)_*W(K(X)) \xrightarrow{V^n} W(K(X)) \to W_n(K(X)) \to 0.$$
In particular, $V^n: (F^n_X)_*(WI_E(p^nD) ) \to WI_E(D) $ is injective. 
By construction, 
the latter homomorphism $WI_E(D)  \to W_nI_E(D) $ is surjective. 
Let us prove the exactness on the middle term $WI_E(D)$. 
Take an element
$$\varphi =(\varphi_0, \varphi_1, \cdots) \in \Gamma(X, WI_E(D))$$
whose image to $W_nI_E(D)$ is zero, 
i.e., $\varphi_0=\cdots=\varphi_{n-1}=0$. 
We can check directly from Definition \ref{d-div-def} that 
$$(\varphi_n, \varphi_{n+1}, \cdots) \in \Gamma(X, WI_E(p^nD)),$$
and hence the sequence (\ref{e-div-induction1}) is exact. 
It holds by the same argument that also (\ref{e-div-induction2}) is exact. 
\end{proof}

\begin{prop}\label{p-div-coherent}
Let $X$ be an integral normal $F$-finite noetherian $\F_p$-scheme. 
Let $D$ be an $\R$-divisor on $X$ and let $E$ be an effective $\R$-divisor. 
Then, for any positive integer $n$, 
$W_nI_E(D)$ is a coherent $W_n\MO_X$-module. 
\end{prop}

\begin{proof}
By Lemma \ref{l-div-sub}(2), 
$WI_E(D)$ is a quasi-coherent $W_n\MO_X$-module. 
Since $X$ is $F$-finite, $(F^n_X)_*(M)$ is a coherent $W_n\MO_X$-module 
for any coherent $W_n\MO_X$-module $M$. 
Thus the assertion follows from induction on $n$ 
and the exact sequence (\ref{e-div-induction2}) of Lemma \ref{p-div-induction}. 
\end{proof}

\begin{prop}
Let $X$ be an integral regular $F$-finite noetherian $\F_p$-scheme. 
Let $D$ be an $\R$-divisor on $X$ and let $E$ be an effective $\R$-divisor. 
Then, for any positive integer $n$, 
$W_nI_E(D)$ is a Cohen-Macaulay $W_n\MO_X$-module. 
\end{prop}

\begin{proof}
Applying the local cohomology to  
\[
0 \to F_*^n(I_E(p^nD)) \xrightarrow{V^n}   W_{n+1}I_E(D) \to W_nI_E(D) \to 0, 
\]
the problem is reduced to the case when $n=1$. 
Then the assertion holds, since $W_1I_E(D) = \MO_X(D-E)$ is an invertible sheaf. 
\end{proof}

\subsection{Weak quasi-$F^e$-splitting}\label{s-weak-QFS}

In this subsection, we define weak quasi-$F^e$-splitting, which generalizes the notion of quasi-$F^e$-splitting in \cite{TWY24}.

\begin{dfn}\label{d B^e}
Let $X$ be an integral normal $F$-finite noetherian $\F_p$-scheme. 
Take an $\R$-divisor $\Delta$ and effective $\R$-divisor $E$. 
Fix positive integers $e$ and $n$.
We define $B^{E, e}_{X, \Delta, n}$ by 
\[
0 \to W_nI_{E}(\Delta) \xrightarrow{F^e} F_*^e W_nI_E(p^e\Delta) \to B^{E, e}_{X, \Delta, n} \to 0.
\]
Thus $ B^{E, e}_{X, \Delta, n}$ is a coherent $W_n\MO_X$-module. 
\end{dfn}

\begin{prop}\label{p Bn ind}
Let $X$ be an integral normal $F$-finite noetherian $\F_p$-scheme. 
Take an $\R$-divisor $\Delta$ and effective $\R$-divisor $E$. 
Then we have the following exact sequence: 
\[
0 \to F_*^n B^{E, e}_{X, p^n\Delta, 1} \to B^{E, e}_{X, \Delta, n+1} \to B^{E, e}_{X, \Delta, n} \to 0
\]
\end{prop}

\begin{proof}
Applying the snake lemma to the following diagram, we obtain the desired exact sequence:  
{\tiny 
\[
\begin{CD}
0 @>>> F_*^{n}(I_{E}(p^{n}\Delta)) @>V^n>>W_{n+1}I_{E}(\Delta) @>R>> W_nI_{E}(\Delta)@>>> 0\\
@. @VVF^eV @VV{F^e}V @VV{F^e}V\\
0 @>>> F_*^{n+e}(I_E(p^{n+e}D)) @>V^n>> F_*^e(W_{n+1}I_E(p^e\Delta)) @>R>> F_*^e(W_{n}I_E(p^e\Delta)) @>>> 0.
\end{CD}
\]
}
\end{proof}

Let $X$ be an integral normal $F$-finite noetherian $\F_p$-scheme. 
Take an $\R$-divisor $\Delta$ and effective $\R$-divisor $E$. 
We define a coherent 
$W_n\MO_X$-module $Q^{E, e}_{X, \Delta, n}$ and a $W_n\MO_X$-module homomorphism $\Phi^{E,e}_{X, \Delta, n}$ by the following pushout diagram: 
\[
\begin{tikzcd}
W_nI_E(\Delta) \arrow[r, "F^e"] \arrow[d, "R^{n-1}"'] & F_*^eW_nI_E(p^e\Delta) \arrow[d]\\
I_E(\Delta) \arrow[r, "\Phi^{E,e}_{X, \Delta, n}"] & Q^{E, e}_{X, \Delta, n}. 
\end{tikzcd}
\]
By definition, we have the following commutative diagram in which all horizontal and vertical sequences are exact: 
\begin{equation}\label{e-big-BQ-diagram}
\begin{tikzcd}
& 0 \arrow{d} & 0 \arrow{d} & & \\
&F_*W_{n-1}I_E(p\Delta) \arrow{d}{V} \arrow[r,dash,shift left=.1em] \arrow[r,dash,shift right=.1em] & F_*W_{n-1}I_E(p\Delta) \arrow{d}{F^eV =VF^e} & & \\
0 \arrow{r} & W_nI_E(\Delta) \arrow{r}{F^e} \arrow{d}{R^{n-1}} & F^e_*W_nI_E(p^e\Delta) \arrow{r} \arrow{d} & B^{E, e}_{X, \Delta, n} \arrow[d,dash,shift left=.1em] \arrow[d,dash,shift right=.1em] \arrow{r}  & 0 \\
0 \arrow{r} & I_E(\Delta) \arrow{r}{\Phi^e_{X, \Delta, n}} \arrow{d} &  \arrow[lu, phantom, "\usebox\pushoutdr" , very near start, yshift=0em, xshift=0.6em, color=black] Q^{E, e}_{X, \Delta, n} \arrow{r} \arrow{d} & B^{E, e}_{X, \Delta, n} \arrow{r} & 0.\\
& 0 & 0 & &
\end{tikzcd}
\end{equation}

\begin{lem}\label{l-Qn-induction}
Let $X$ be an integral normal $F$-finite noetherian $\F_p$-scheme. 
Take an $\R$-divisor $\Delta$ and effective $\R$-divisor $E$. 
Then we have an exact sequence
\[
0 \to F_*^n B^{E, e}_{X, p^n\Delta, 1} \to Q^{E, e}_{X, \Delta, n+1} \to Q^{E, e}_{X,\Delta, n} \to 0. 
\]
\end{lem}

\begin{proof}
The assertion holds by applying the snake lemma to the following commutative diagram in which each horizontal sequence is exact: 
\[
\begin{tikzcd}
0 \arrow{r} & F^n_*I_E(p^n\Delta) \arrow{r}{V^{n-1}} \arrow{d}{F^e} & 
F_*W_{n}I_E(p\Delta)\arrow[r, "R"] \arrow{d}{F^eV=VF^e} & 
F_*W_{n-1}I_E(p\Delta) \arrow[d, "F^eV=VF^e"] \arrow{r}  & 0 \\
0 \arrow{r} & 
F_*^{e+n}I_E(p^{e+n}\Delta) \arrow{r}{V^n}  &   
F^e_*W_{n+1}I_E(p^e\Delta)  \arrow[r, "R"] & 
F^e_*W_nI_E(p^e\Delta)  \arrow{r} & 0.
\end{tikzcd}
\]
\end{proof}

\begin{prop}\label{p-BQ-CM}
Let $X$ be an integral regular $F$-finite noetherian $\F_p$-scheme. 
Assume that 
$F^e : I_E(\Delta)_x \to F_*^eI_E(p^e\Delta)_x$ splits as an $\MO_{X, x}$-module homomorphism 
for every integer $e>0$ and every point $x \in X$. 
Then 
$B^{E, e}_{X, \Delta, n}$ and $Q^{E, e}_{X, \Delta, n}$ are coherent Cohen-Macaulay 
$W_n\MO_X$-modules. 
\end{prop}

\begin{proof}
By (\ref{e-big-BQ-diagram}), it is enough to show that $B^{E, e}_{X, \Delta, n}$ 
is Cohen-Macaulay. 
By the exact sequence (Proposition \ref{p Bn ind}): 
\[
0 \to F_*^n B^{E, e}_{X, p^n\Delta, 1} \to B^{E, e}_{X, \Delta, n+1} \to B^{E, e}_{X, \Delta, n} \to 0,
\]
we may assume that $n=1$, i.e., it suffices to prove that $B^{E, e}_{X, \Delta, 1}$ is 
Cohen-Macaulay. 
Fix $x \in X$. 
By definition, we have the following exact sequence: 
\[
0 \to I_E(\Delta)_x \xrightarrow{F^e} F_*^eI_E(p^e\Delta)_x \to (B^{E, e}_{X, \Delta, 1})_x \to 0. 
\]
By assumption, this exact sequence splits. Thus, $(B^{E, e}_{X, \Delta, 1})_x $ 
is a direct summand of $F_*^eI_E(p^e\Delta)_x$, which shows that $(B^{E, e}_{X, \Delta, 1})_x$ is 
Cohen-Macaulay.
\end{proof}

For the $W_n\MO_X$-module homomorphism 
\[
\Phi^{E, e}_{X, \Delta, n} : I_E(\Delta) \to Q^{E, e}_{X, \Delta, n}, 
\]
we apply the contravariant functor $(-)^* := \cHom_{W_n\MO_X}(-, W_n\omega_X(-K_X))$: 
\[
(\Phi^{E, e}_{X, \Delta, n})^* : (Q^{E, e}_{X, \Delta, n})^* \to I_E(\Delta)^*, 
\]
where $W_n\omega_X(-K_X)$ denotes the $S_2$-hull of $W_n\omega_X \otimes W_n\MO_X(-K_X)$ (cf.\ \cite{TWY24}). 
Note that 
\[
I_E(\Delta)^* =  \cHom_{W_n\MO_X}(I_E(\Delta), W_n\omega_X(-K_X)) 
\simeq \cHom_{\MO_X}(I_E(\Delta), \MO_X) \simeq \MO_X(-\rdown{\Delta -E}).
\]
In particular, if $\Delta = E$, then $I_E(\Delta)^* \simeq \MO_X$

\begin{dfn}\label{d QFS}
Let $X$ be a normal variety and $\Delta$ an effective $\R$-divisor on $X$.
Fix a positive integer $e>0$.
For a positive integer $n>0$, we say that $(X, \Delta)$ is {\em weakly $n$-quasi-$F^e$-split} if the map
the map
\[
H^0(X, (\Phi^{\Delta, e}_{X, \Delta, n})^*) : 
H^0(X, (Q^{\Delta, e}_{X, \Delta, n})^*) \to H^0(X, I_{\Delta}(\Delta)^*) 
\]
is surjective. 
Recall that $H^0(X, I_{\Delta}(\Delta)^*) \cong H^0(X, \MO_X)$. 
When $\Delta=0$, we simply say that $X$ is {\em $n$-quasi-$F^e$-split}.

We say $(X, \Delta)$ (resp.~$X$) is \textit{weakly quasi-$F^e$-split} (resp.~\textit{quasi-$F^e$-split}) if there exists a positive integer $n>0$ such that
$(X, \Delta)$ (resp.~$X$) is $n$-weakly quasi-$F^e$-split (resp.~$n$-quasi-$F^e$-split).
\end{dfn}

\begin{rem}
    The above definition of quasi-$F^e$-splitting coincides with that in \cite{TWY24}. 
\end{rem}

\subsection{Inversion of adjunction for log Calabi-Yau pairs}\label{ss-logCY-IOA}

The purpose of this subsection is to establish an inversion of adjunction of quasi-$F$-splitting for log Calabi-Yau pairs (Corollary \ref{c-IOA}). 

\begin{notation}\label{n-IOA}
Let $\kappa$ be an $F$-finite field of characteristic $p>0$. 
Let $X$ be a regular projective variety over $\kappa$. 
Take a regular prime divisor $S$ on $X$. 
Let $B$ be an effective $\R$-divisor such that 
$S \not\subset \Supp\,B$, all the coefficients of $B$ are contained in $[0, 1]$, 
and $S+B$ is simple normal crossing around $S$ 
(i.e., for every closed point $x \in S$, 
there exist $1 \leq r \leq \dim X$ and a regular system of parameters $x_1, ..., x_{\dim X}$ of 
the regular local ring $\MO_{X, x}$ such that $\Supp (S+B)|_{\Spec \MO_{X, x}} = \Spec \MO_{X, x} / (x_1 \cdots x_r)$). 
Set $B_S \coloneqq B|_S$. 
\end{notation}

\begin{lem}\label{l-IOA-diag}
We use Notation \ref{n-IOA}. 
Fix positive integers $e>0$ and  $n>0$. 
Assume that 
\begin{enumerate}
    \item $(S, B_S)$ is weakly $n$-quasi-$F^e$-split, and 
    \item $H^1(X, 
    \mathcal Hom_{W_n\MO_X}(Q^{B, e}_{X, B, n} \otimes_{W_n\MO_X} W_n\MO_X(S), W_n\omega_X(-K_X)))=0$. 
\end{enumerate}
Then $(X, S+B)$ is weakly $n$-quasi-$F^e$-split. 
\end{lem}

\begin{proof}
By \cite[Proposition 2.28]{KTTWYY1}, we have the following commutative diagram in which each horizontal sequence is a natural exact sequence: 
\[
\begin{tikzcd}
0 \arrow[r] &  F_*^eW_nI_{X, S+B}(p^eB) \arrow[r] \arrow[d, leftarrow, "F^e"'] & F_*^eW_nI_{X, B}(p^eB) \arrow[r] \arrow[d, leftarrow, "F^e"'] &  F_*^eW_nI_{S, B_S}(p^eB_S) \arrow[r] \arrow[d, leftarrow, "F^e"'] & 0 \\
0 \arrow[r] & W_nI_{X, S+B}(B) \arrow[r] & W_nI_{X, B}(B) \arrow[r] & W_nI_{S, B_S}(B_S) \arrow[r] & 0. 
\end{tikzcd}
\]
Note that the left vertical arrow is obtained by applying $(-) \otimes W_n\MO_X(-S)$ to 
\[
F^e : W_nI_{X, S+B}(S+ B) \to 
F_*^eW_nI_{X, S+B}(p^e(S+B)). 
\]
Taking pushouts, we obtain another exact sequence in which each horizontal sequence is exact: 
\[
\begin{tikzcd}
0 \arrow[r] &  
Q^{S+B, e}_{X, S+B, n} \otimes W_n\MO_X(-S)
\arrow[r] \arrow[d, leftarrow, "\Phi'"'] &
Q^{B, e}_{X, B, n} \arrow[r] \arrow[d, leftarrow, "\Phi^{B, e}_{X, B, n}"'] &  
Q^{B_S, e}_{S, B_S, n} \arrow[r] \arrow[d, leftarrow, "\Phi^{B_S, e}_{S, B_S, n}"'] & 0 \\
0 \arrow[r] & I_{X, S} \arrow[r] & \MO_X \arrow[r] & \MO_{S} \arrow[r] & 0. 
\end{tikzcd}
\]
We now apply  $\mathcal Hom(-, W_n\omega_X) \otimes W_n\MO_X(-(K_X+S))$, 
which is the same as $( - \otimes_{W_n\MO_X} W_n\MO_X(S))^*$ for $(-)^* := \mathcal Hom_{W_n\MO_X}(-, W_n\omega_X(-K_X))$. 
We then  obtain an  exact sequence in which each horizontal sequence is exact: 
\[
\begin{tikzcd}
0 \arrow[r] &  
(Q^{B, e}_{X, B, n} \otimes_{W_n\MO_X} W_n\MO_X(S))^* \arrow[r] \arrow[d, "(\Phi^{B, e}_{X, B, n} \otimes W_n\MO_X(S))^*"'] &
(Q^{S+B, e}_{X, S+B, n})^* \arrow[r] \arrow[d, "(\Phi^{S+B, e}_{X, S+B, n})^*"'] &  
(Q^{B_S, e}_{S, B_S, n})^* \arrow[r] \arrow[d, "(\Phi^{B_S, e}_{S, B_S, n})^*"'] & 0 \\
0 \arrow[r] & I_{X, S} \arrow[r] & \MO_X \arrow[r] & \MO_S \arrow[r] & 0. 
\end{tikzcd}
\]
Here we used 
\begin{itemize}
\item $(I_{X, S} \otimes_{W_n\MO_X} W_n\MO_X(S))^*= \MO_X^* \simeq \MO_X$, 
\item $(\MO_X \otimes_{W_n\MO_X} W_n\MO_X(S))^*= (\MO_X(S))^* \simeq \MO_X(-S)$, and 
\item $(\MO_S \otimes_{W_n\MO_X} W_n\MO_X(S))^* \simeq \MO_S$, 
where the last isomorphism follows from 
\[
(\MO_S \otimes_{W_n\MO_X} W_n\MO_X(S))^*
\simeq \mathcal Hom_{W_n\MO_X}(\MO_S, W_n\omega_X) \otimes_{W_n\MO_X} W_n(-(K_X+S))
\]
\[
\simeq \mathcal Hom_{\MO_S}(\MO_S, \omega_S) \otimes_{\MO_S} \MO_S(-K_S) \simeq \MO_S. 
\]
\end{itemize}
Taking cohomologies, we obtain the following commutative diagram in which the upper sequence is exact: 
\[
\begin{tikzcd}
H^0(X, (Q^{S+B, e}_{X, S+B, n})^*) \arrow[r, "\alpha"] \arrow[d, "{H^0(X, (\Phi^{S+B, e}_{X, S+B, n})^*)}"'] &  
H^0(S, (Q^{B_S, e}_{S, B_S, n})^*) \arrow[r] \arrow[d, "{H^0(S, (\Phi^{B_S, e}_{S, B_S, n})^*)}"'] & 
H^1(X, (Q^{B, e}_{X, B, n} \otimes_{W_n\MO_X} W_n\MO_X(S))^*)
\overset{{\rm (2)}}{=}0 
\\
H^0(X,  \MO_X) \arrow[r, "\beta"] & H^0(S, \MO_S). 
\end{tikzcd}
\]
By (1), $H^0(S, (\Phi^{B_S, e}_{S, B_S, n})^*)$ is surjective 
(Definition \ref{d QFS}). 
Then the composite arrow of the above square diagram is surjective. 
Then $\beta$ is a surjective ring homomorphism from a field $H^0(X, \MO_X)$, 
and hence $\beta$ is bijective. 
Therefore, $H^0(X, (\Phi^{S+B, e}_{X, S+B, n})^*)$ is surjective, 
i.e., $(X, S+B)$ is weakly $n$-quasi-$F^e$-split (Definition \ref{d QFS}). 
\end{proof}

\begin{thm}\label{t-IOA}
We use Notation \ref{n-IOA}. 
Fix integers $e>0$ and $n>0$.  
Assume that the following hold. 
\begin{enumerate}
\item $p^e(K_X+S+B) \sim 0$. 
\item $H^0(X, \MO_X(-S + \ulcorner  (1-p^m)(K_X+S+B) \urcorner))=0$ for every integer $m >0$. 
\item $H^1(X, \MO_X(K_X+ \ulcorner B \urcorner))=0$. 
\end{enumerate} 
If $(S, B_S)$ is weakly $n$-quasi-$F^e$-split, then so is $(X, S+B)$. 
\end{thm}

\begin{proof}
Assume that $(S, B_S)$ is weakly $n$-quasi-$F^e$-split. 
By Lemma \ref{l-IOA-diag}, it suffices to show 
\[
H^1(X, 
    \mathcal Hom_{W_m\MO_X}(Q^{B, e}_{X, B, m} \otimes_{W_m\MO_X} W_m\MO_X(S), W_m\omega_X(-K_X)))=0
\] 
for every $1 \leq m \leq n$. 
Let us prove this by induction on $m$. 
If $m=1$, then 
\begin{eqnarray*}
&&H^1(X, 
    \mathcal Hom_{W_m\MO_X}(Q^{B, e}_{X, B, m} \otimes_{W_m\MO_X} W_m\MO_X(S), W_m\omega_X(-K_X)))\\
&=& H^1(X, 
    \mathcal Hom_{\MO_X}(F_*^eI_B(p^eB)\otimes_{\MO_X} \MO_X(S), \omega_X(-K_X)))\\
&\simeq & H^1(X, \mathcal Hom_{\MO_X}(F_*^e(I_B \otimes_{\MO_X} \MO_X(p^e(K_X+S+B)), \omega_X)\\
&\simeq & H^1(X, \mathcal Hom_{\MO_X}(F_*^eI_B, \omega_X)\\
&\simeq & H^1(X, F_*^e\mathcal Hom_{\MO_X}(\MO_X(-\rup{B}), \omega_X)\\
&\simeq & H^1(X, \MO_X(K_X + \rup{B}))\overset{(3)}{=}0. 
\end{eqnarray*}
Fix an integer $m$ satisfying $1 \leq m <n$. 
By the induction hypothesis and the exact sequence:
\[
0 \to F_*^m B^{B, e}_{X, p^mB, 1} \to Q^{B, e}_{X, B, m+1} \to Q^{B, e}_{X, B, m} \to 0, 
\]
it suffices to show that 
\[
H^1(X, \mathcal Hom(F_*^m B^{B, e}_{X, p^mB, 1}\otimes_{W_m\MO_X} W_m\MO_X(S), W_m\omega_X(-K_X)))=0. 
\]
By Definition \ref{d B^e}, we have 
\begin{equation}\label{e1-IOA}
0 \to F_*^mI_B(p^mB) \xrightarrow{F^e} F_*^{m+e}(I_B(p^{m+e}B)) \to F_*^mB^{B, e}_{X, p^mB, 1} \to 0. 
\end{equation}
In particular, $F_*^m B^{B, e}_{X, p^mB, 1}$ is an $\MO_X$-module, and hence 
\begin{eqnarray*}
&& \mathcal Hom_{W_m\MO_X}(F_*^m B^{B, e}_{X, p^mB, 1}\otimes_{W_m\MO_X} W_m\MO_X(S), W_m\omega_X(-K_X))\\
&=& \mathcal Hom_{W_m\MO_X}(F_*^m B^{B, e}_{X, p^mB, 1}\otimes_{\MO_X} \MO_X(K_X+S), W_m\omega_X)\\
&\simeq & \mathcal Hom_{\MO_X}(F_*^m B^{B, e}_{X, p^mB, 1}\otimes_{\MO_X} \MO_X(K_X+S), \omega_X). 
\end{eqnarray*}
Thus the problem is reduced to showing that 
\[
H^1(X, \mathcal Hom_{\MO_X}(F_*^m B^{B, e}_{X, p^mB, 1}\otimes_{\MO_X} \MO_X(K_X+S), \omega_X)=0.
\]
Applying $\mathcal Hom_{\MO_X}(- \otimes_{\MO_X} \MO_X(K_X+S), \omega_X)$ 
to (\ref{e1-IOA}), we obtain  an exact sequence
\begin{eqnarray*}
0 
&\to& \cHom_{\MO_X}((F_*^mB^{B, e}_{X, p^mB, 1} \otimes \MO_X(K_X+S), \omega_X)\\
&\to& \cHom_{\MO_X}((F_*^{m+e}(I_B(p^{m+e}B)) \otimes \MO_X(K_X+S), \omega_X)\\
&\to& \cHom_{\MO_X}(F_*^mI_B(p^mB) \otimes \MO_X(K_X+S), \omega_X) \to 0.
\end{eqnarray*}
Therefore, it is enough to show that 
\begin{enumerate}
\item[(I)] $H^0(X, \cHom_{\MO_X}(F_*^mI_B(p^mB) \otimes \MO_X(K_X+S), \omega_X))=0$, and 
\item[(II)] $H^1(X, \cHom_{\MO_X}((F_*^{m+e}(I_B(p^{m+e}B)) \otimes \MO_X(K_X+S), \omega_X))=0$. 
\end{enumerate}

Let us show (I). 
We have  that 
\begin{eqnarray*}
&&H^0(X, \cHom_{\MO_X}(F_*^mI_B(p^mB) \otimes \MO_X(K_X+S), \omega_X)) \\
&\simeq &H^0(X, \cHom_{\MO_X}(F_*^mI_B(p^m(K_X+S+B)), \omega_X)) \\
&\simeq&H^0(X, F_*^m \cHom_{\MO_X}(\MO_X(p^m(K_X+S+B)-B), \omega_X))\\
&\simeq& H^0(X, \cHom_{\MO_X}(\MO_X(p^m(K_X+S+B)-(K_X+B)), \MO_X)\\
&=& H^0(X, \cHom_{\MO_X}(\MO_X(S+(p^m-1)(K_X+S+B), \MO_X)\\
&\simeq& H^0(X, \MO_X( -S +\rup{(1-p^m)(K_X+S+B)})\overset{(2)}{=}0. 
\end{eqnarray*}
Thus (I) holds.

Let us show (II). 
By the same argument as in (I), we obtain 
\begin{eqnarray*}
&&H^1(X, \cHom_{\MO_X}((F_*^{m+e}(I_B(p^{m+e}B)) \otimes \MO_X(K_X+S), \omega_X)) \\
&\simeq& H^1(X, \cHom_{\MO_X}(\MO_X(p^{m+e}(K_X+S+B)-(K_X+B)), \MO_X)\\
&\overset{{\rm (1)}}{\simeq}& H^1(X, \cHom_{\MO_X}(\MO_X(-(K_X+B)), \MO_X)\\
&=& H^1(X, \MO_X( K_X + \rup{B})\overset{(3)}{=}0. 
\end{eqnarray*}
Thus (II) holds, and we conclude that $(X,S+B)$ is weakly $n$-quasi-$F^e$-split.
\end{proof}

\begin{cor}\label{c-IOA}
Let $\kappa$ be an $F$-finite field of characteristic $p>0$. 
Let $X$ be a regular projective variety over $\kappa$. 
Take a regular prime divisor $S$ on $X$. 
Let $B$ be a reduced divisor such that 
$S \not\subset \Supp\,B$ and $S+B$ is simple normal crossing around $S$.  
Set $B_S \coloneqq B|_S$. 
Assume that 
\begin{enumerate}
\item $K_X+S+B \sim 0$. 
\item $H^1(X, \MO_X(-S))=0$. 
\end{enumerate}
Fix integers $e>0$ and $n>0$.  
If $(S, B_S)$ is weakly $n$-quasi-$F^e$-split, then so is $(X, S+B)$. 
\end{cor}

\begin{proof}
It is enough to check that the assumptions (1)-(3) in Theorem \ref{t-IOA} hold. 
Clearly, (1) implies Theorem \ref{t-IOA}(1), and (1) and (2) imply Theorem \ref{t-IOA}(2). 
Theorem \ref{t-IOA}(3) follows from (2), because 
$K_X + \rup{B} = K_X +B \sim -S$. 
\end{proof}

%% file: section6.tex
\section{Quasi-$F$-splitting and global $F$-regularity of smooth del Pezzo varieties}\label{s-GFR}

\subsection{Quasi-$F$-splitting}\label{ss-pf-main-thm}

In this subsection, we prove that smooth del Pezzo varieties are quasi-$F$-split, 
which is the main result of this article. 

\begin{prop}\label{p-IOA-Fano3}
Let $k$ be an algebraically closed field of characteristic $p>0$ and let $k \subset \kappa$ be a field extension, where $\kappa$ is $F$-finite. 
Let $X$ be an $n$-dimensional  smooth projective variety over $k$. 
Assume that there exist prime divisors $H_1,\ldots,H_{n-1}$ on $X \times_k \kappa$ such that the following hold. 
\begin{enumerate}
\renewcommand{\labelenumi}{(\alph{enumi})}
\item $K_{X \times_k \kappa}+H_1+\cdots+ H_{n-1} \sim 0$. 
\item 
$H_{j_1} \cap \cdots \cap H_{j_i}$ is regular and integral 
if $1 \leq i \leq n-1$ and $j_1 < \cdots < j_i$. 
\item 
$H^i(X \times_k \kappa, \MO_{X \times_k \kappa}(-H_{j_1}- \cdots -H_{j_{i}}))=0$ if 
$1 \leq i \leq n-1$ and $j_1 < \cdots < j_i$. 
\item $\bigcap_{i=1}^{n-1} H_i$ is one-dimensional and smooth over $\kappa$. 
\end{enumerate}
Then the following hold. 
\begin{enumerate}
    \item $(X \times_k \kappa, H_1+\cdots+ H_{n-1})$ is weakly $2$-quasi-$F$-split. 
    \item $H_{j_1} \cap \cdots \cap H_{j_i}$ is $2$-quasi-$F$-split 
if $1 \leq i \leq n-1$ and $j_1 < \cdots < j_i$. 
\item $X$ is $2$-quasi-$F$-split. 
\end{enumerate}
\end{prop}

\begin{proof}
Set 
\[
Y_n := X \times_k \kappa,\quad 
Y_{n-1} := H_{n-1}, \quad 
Y_{n-2} := H_{n-2} \cap H_{n-1}, ..., Y_1 := H_1 \cap \cdots \cap H_{n-1}. 
\]
For each $i$, we define an effective Cartier divisor $B_j$ on $Y_j$ as follows: 
\[
B_n := H_1 + \cdots + H_{n-2}, \quad 
B_{n-1} :=  (H_1 + \cdots + H_{n-3})|_{Y_{n-1}}
\]
\[
B_{n-2} := (H_1 + \cdots + H_{n-4})|_{Y_{n-2}},...,  
D_{3} :=  H_1|_{Y_{3}}, \quad  D_2 := 0, \quad D_1 :=0. 
\]
Clearly, it holds that 
\[
(K_{Y_j} + Y_{j-1} + B_{j})|_{Y_{j-1}} = K_{Y_{j-1}} +(B_j|_{Y_{j-1}}) 
= K_{Y_{j-1}}  +Y_{j-2} + B_{j-1}. 
\]
This, together with (a), implies  $K_{Y_j} + Y_{j-1} + B_{j} =(K_{Y_n} + Y_{n-1} +B_n)|_{Y_j} \sim 0$.

\begin{claim*}
$H^1(Y_j, \MO_{Y_j}(-Y_{j-1}))=0$ for every $2 \leq j\leq n$. 
\end{claim*}

\begin{proof}[Proof of Claim]
Fix  $2 \leq j\leq n$. 
Recall that $Y_{j-1} = {Y_j} \cap H_{j-1} = H_{j-1}|_{Y_j}$. 
By the exact sequence
\[
0 \to \MO_{Y_{j+1}}(-H_{j-1} -H_j) \to \MO_{Y_{j+1}}(-H_{j-1}) \to \MO_{Y_j}(-H_{j-1})\to 0, 
\]
it is enough to show 
\[
H^2(Y_{j+1}, \MO_{Y_{j+1}}(-H_{j-1} -H_j)) =H^1(Y_{j+1}, \MO_{Y_{j+1}}(-H_{j-1}))=0. 
\]
By the exact sequence
\[
0 \to \MO_{Y_{j+2}}(-H_{j+1}) \to \MO_{Y_{j+2}} \to \MO_{Y_{j+1}} \to 0, 
\]
it suffices to prove 
\[
H^2(Y_{j+2}, \MO_{Y_{j+2}}(-H_{j-1} -H_j)) =
H^3(Y_{j+2}, \MO_{Y_{j+2}}(-H_{j-1} -H_j-H_{j+1})) =
\]
\[
=H^1(Y_{j+2}, \MO_{Y_{j+2}}(-H_{j-1})) =
H^2(Y_{j+2}, \MO_{Y_{j+2}}(-H_{j-1} -H_{j+1})) =0. 
\]
Repeating this procedure, the problem is reduced to showing that 
\[
H^i(Y_n, \MO_{Y_n}(-H_{j_1} - \cdots -H_{j_i})) =0 
\]
for  $1 \leq i \leq n-j+1$ and $j-1 =j_1 <j_2 < \cdots < j_i \leq n-1$. 
By $2 \leq j \leq n$, we have $n-j +1 \leq n-1$ and $1 \leq j-1$. 
Hence $1 \leq i \leq n-1$ and $1 \leq j_1 < j_2 < \cdots < j_i \leq n-1$. 
Therefore, the required vanishing is assured by (c). 
This completes the proof of Claim. 
\end{proof}

Set $Y_0 \coloneqq 0$. 
It is enough to prove $(\star)_j$ for every $1 \leq j \leq n$. 
\begin{enumerate}
    \item[$(\star)_j$] $(Y_j, Y_{j-1} +B_j)$ is weakly $2$-quasi-$F$-split. 
\end{enumerate}
Indeed, $(\star)_n$ implies (1), (2) holds by symmetry, and (3) follows from (1) \cite[Proposition 2.12]{KTY}. 
We prove $(\star)_j$ by induction on $j$.

Let us show $(\star)_1$. 
Since $Y_1 = \bigcap_{i=1}^{n-1} H_i$ is a smooth projective one-dimensional integral scheme over $\kappa$ with $K_{Y_1} \sim 0$, 
the base change $Y_1 \times_{\kappa} \overline \kappa$ 
to the algebraic closure $\overline \kappa$ of $\kappa$ is a disjoint union of elliptic curves, which is $2$-quasi-$F$-split (\cite[Remark 2.11]{KTTWYY1}).   
Thus so is $Y_1$ \cite[Proposition 2.12]{KTY}. 
This completes the proof of $(\star)_1$.


Fix $1 < i \leq n$ and assume $(\star)_{j-1}$, i.e., 
$(Y_{j-1}, Y_{j-2} +B_{j-1})$ is weakly $2$-quasi-$F$-split. 
Since $K_{Y_{j}} + Y_{j-1} + B_{j} \sim 0$ and $H^1(Y_{j}, \MO_{Y_{j}}(-Y_{j-1}))=0$ by Claim, 
it follows from Corollary \ref{c-IOA} that 
$(\star)_{j-1}$ implies $(\star)_j$. 
Thus $(\star)_j$ holds. 
\end{proof}

\begin{thm}\label{mainthm:quasi-F-split}
    Let $(X,L)$ be a smooth del Pezzo variety over an algebraically closed field $k$ of characteristic $p>0$.
    Then $X$ is $2$-quasi-$F$-split.
\end{thm}

\begin{proof}
Set $n\coloneqq\dim X$. 
For $m\coloneqq  n-1, k^{(0)} \coloneqq k, X^{(0)} \coloneqq X, L^{(0)}\coloneqq L$, 
    we use the same notation as in (\ref{n generic repeated}). 
    It is enough to  confirm that $X$ and $H_{1},\ldots, H_{n-1}$ satisfy the conditions (a)-(d) of Proposition \ref{p-IOA-Fano3}.
    The conditions (a) and (b) of Proposition \ref{p-IOA-Fano3} follows from  
    (\ref{n generic repeated})(iv) and Lemma \ref{l generic ladder}. 
    By the definition of del Pezzo varieties (Definition \ref{d dP var}), we have
    \[
    H^i(X \times_k \kappa, \sO_{X \times_k \kappa}(-H_{j_1} - \cdots -H_{j_i}))=H^i(X , \sO_{X}(-iL))\otimes_{k}\kappa =0
    \]
    for all $i\in\{1,\ldots,n-1\}$ and $j_1 < \cdots < j_i$, and thus (c) of Proposition \ref{p-IOA-Fano3} is satisfied. 
    
    Finally, let us show (d) of Proposition \ref{p-IOA-Fano3}. 
    Since $\bigcap_{i=1}^{n-1}H_{i}$ is one-dimensional 
    (Lemma \ref{l generic ladder}), it is enough to show that 
    $\bigcap_{i=1}^{n-1}H_{i}$ is smooth. 
    By Proposition \ref{p generic vs general}, 
    it is enough to find $D_1,\ldots, D_{n-1} \in |L|$ such that 
    $D_1 \cap \cdots \cap D_{n-1}$ is smooth. 
    This follows from Theorem \ref{thm:normality of ladder}. 
\end{proof}

\subsection{Global $F$-regularity}\label{subsection:6-2}

In this subsection, we focus on the global $F$-regularity of del Pezzo varieties.

\begin{dfn}
Let $X$ be a normal variety and let $\Delta$ be an effective $\Q$-divisor on $X$. 
\begin{enumerate}
\item We say that $(X, \Delta)$ is {\em $F$-split} if 
\[
\MO_X \xrightarrow{F^e} F_*^e\MO_X \hookrightarrow F_*^e\MO_X( \rdown{(p^e-1)\Delta})
\]
splits as an $\MO_X$-module homomorphism for every $e \in \Z_{>0}$. 
\item 
We say that $(X, \Delta)$ is {\em globally $F$-regular} 
if, given an effective $\Z$-divisor $E$, there exists $e \in \Z_{>0}$ such that 
\[
\MO_X \xrightarrow{F^e} F_*^e\MO_X \hookrightarrow 
F_*^e\MO_X( \rup{(p^e-1)\Delta} +E)
\]
splits as an $\MO_X$-module homomorphism. 
\item 
We say that $(X, \Delta)$ is {\em strongly $F$-regular} 
if, for every closed point $x\in X$, there exists an open subscheme $U$ of $X$ 
such that $x \in U$ and $(U, \Delta|_U)$ is globally $F$-regular. 
\end{enumerate}
We say that $X$ is {\em $F$-split} (resp. {\em globally $F$-regular}) 
if so is $(X, 0)$. 
\end{dfn}

\begin{rem}\label{rem:gfr}\,
\begin{enumerate}
    \item By definition, globally $F$-regular varieties are $F$-split.
When $X$ is a normal projective strongly $F$-regular variety such that $-K_X$ is $\Q$-Cartier and ample, then the converse direction holds (the proof of \cite[Theorem 6.2]{KT23}).
\item By definition, globally $F$-regular varieties are strongly $F$-regular.
      It is known that $\Q$-Gorenstein strongly $F$-regular varieties are klt \cite[Theorem 3.3]{Hara-Watanabe}.
\end{enumerate}
\end{rem}

Now, we prove Theorem 
\ref{Intro:F-split}.

\begin{proof}[Proof of Theorem \ref{Intro:F-split}]
By induction on $i$, we prove that $X_i$ is globally $F$-regular for every $i \geq 2$. 
We have the following four cases. 
    \begin{enumerate}
        \item[(i)] $L^n\geq 4$. 
        \item[(ii)] $L^n=3$ and $p>2$. 
        \item[(iii)] $L^n=2$ and $p>3$. 
        \item[(iv)] $L^n=1$ and $p>5$.
    \end{enumerate}
    We first show that $X_2$ is globally $F$-regular in each case.
    Note that $L_2 \sim -K_{X_2}$.
    If (i), (ii), or (iii) holds, $X_2$ is smooth by Proposition \ref{prop:degree geq 3}.
    In this case, $X_2$ is globally $F$-regular by \cite[Example 5.5]{Hara(two-dim)}.
    If (iv) holds, $X_2$ is a canonical del Pezzo surface by Theorem \ref{thm:X_2 is a can. dP}, and thus 
    it is globally $F$-regular by \cite[Theorem A]{Kawakami-Tanaka(dPsurface)}.

Fix  $i \geq 3$ and assume that $X_{i-1}$ is globally $F$-regular. 
It is enough to show that $X_i$ is globally $F$-regular. 
Take an ample effective Cartier divisor $D_i$ on $X_i$ such that $X_i\setminus D_i$ is smooth and affine. 
        Since $X_{i-1}\in |L_i|$ is general, we may assume that $X_{i-1}\not\subset \Supp(D_{i})$.
        Set $D_{i-1}\coloneqq D_i|_{X_{i-1}}$. 
        We can find an integer $e\gg 0$ such that  (a) and (b) below hold for $\epsilon\coloneqq 1/(p^e-1)$. 
        \begin{enumerate}
            \item[(a)] $(X_{i-1},\epsilon D_{i-1})$ is globally $F$-regular (\cite[Corollary 6.1]{SS10}).
            \item[(b)] $-(K_{X_i}+X_{i-1}+\epsilon D_i)$ is ample.
        \end{enumerate}
        Since $(p^e-1)(K_{X_i}+X_{i-1}+\epsilon D_i)$ is Cartier and $(K_{X_i}+X_{i-1}+\epsilon D_i)|_{X_{i-1}}=K_{X_{i-1}}+\epsilon D_{i-1}$, we can apply \cite[Lemma 2.7]{CTW17} to conclude that
        $(X_i, X_{i-1}+\epsilon D_i)$ is $F$-split. 
        In particular, $(X_i, \epsilon D_i)$ is $F$-split, 
        which implies that 
        \[
        \sO_{X_i} \to F_{*}^e\sO_{X_i}((p^e-1)\epsilon D_i)=F_{*}^e\sO_{X_i}(D_i)
        \]
        splits.
        Since $X_i\setminus D_i$ is smooth and affine,
        we conclude that $X_i$ is globally $F$-regular by \cite[Theorem 3.9]{SS10}.
\end{proof}

\begin{cor}\label{cor:ladder is canonical}
    Let $(X,L)$ be a smooth del Pezzo variety.
    Let
    \[
   (X,L)=(X_n,L_n) \supset (X_{n-1},L_{n-1}) \supset \cdots \supset (X_1,L_1)
    \] 
    be a general ladder.
    If $p>5$, then $X_i$ is canonical and strongly $F$-regular for all $i\geq 1$.
\end{cor}
\begin{proof}
    Since $X_1$ is smooth by Theorem \ref{introthm:rungs}, we may assume that $i\geq 2$.
    By Theorem \ref{Intro:F-split}, $X_i$ is globally $F$-regular. 
    Then $X_i$ is klt (Remark \ref{rem:gfr}(2)). 
    Since $X_i$ is Gorenstein, $X_i$ is canonical. 
\end{proof}

%% file: main.bbl

%% file: main.bbl
\begin{bibdiv}
\begin{biblist}

\bib{BT22}{article}{
      author={Bernasconi, Fabio},
      author={Tanaka, Hiromu},
       title={On del {P}ezzo fibrations in positive characteristic},
        date={2022},
        ISSN={1474-7480},
     journal={J. Inst. Math. Jussieu},
      volume={21},
      number={1},
       pages={197\ndash 239},
         url={https://doi-org.utokyo.idm.oclc.org/10.1017/S1474748020000067},
      review={\MR{4366337}},
}

\bib{CD89}{book}{
      author={Cossec, Fran\c{c}ois~R.},
      author={Dolgachev, Igor~V.},
       title={Enriques surfaces. {I}},
      series={Progress in Mathematics},
   publisher={Birkh\"{a}user Boston, Inc., Boston, MA},
        date={1989},
      volume={76},
        ISBN={0-8176-3417-7},
         url={https://doi-org.utokyo.idm.oclc.org/10.1007/978-1-4612-3696-2},
      review={\MR{986969}},
}

\bib{CLS11}{book}{
      author={Cox, David~A.},
      author={Little, John~B.},
      author={Schenck, Henry~K.},
       title={Toric varieties},
      series={Graduate Studies in Mathematics},
   publisher={American Mathematical Society, Providence, RI},
        date={2011},
      volume={124},
        ISBN={978-0-8218-4819-7},
         url={https://doi.org/10.1090/gsm/124},
      review={\MR{2810322}},
}

\bib{CTW17}{article}{
      author={Cascini, Paolo},
      author={Tanaka, Hiromu},
      author={Witaszek, Jakub},
       title={On log del {P}ezzo surfaces in large characteristic},
        date={2017},
        ISSN={0010-437X,1570-5846},
     journal={Compos. Math.},
      volume={153},
      number={4},
       pages={820\ndash 850},
         url={https://doi.org/10.1112/S0010437X16008265},
      review={\MR{3621617}},
}

\bib{Eke87}{incollection}{
      author={Ekedahl, Torsten},
       title={Foliations and inseparable morphisms},
        date={1987},
   booktitle={Algebraic geometry, {B}owdoin, 1985 ({B}runswick, {M}aine, 1985)},
      series={Proc. Sympos. Pure Math.},
      volume={46, Part 2},
   publisher={Amer. Math. Soc., Providence, RI},
       pages={139\ndash 149},
         url={https://doi.org/10.1090/pspum/046.2/927978},
      review={\MR{927978}},
}

\bib{Fujita75}{article}{
      author={Fujita, Takao},
       title={On the structure of polarized varieties with {$\Delta $}-genera zero},
        date={1975},
        ISSN={0040-8980},
     journal={J. Fac. Sci. Univ. Tokyo Sect. IA Math.},
      volume={22},
       pages={103\ndash 115},
      review={\MR{369363}},
}

\bib{Fujita(ch=0)}{article}{
      author={Fujita, Takao},
       title={On the structure of polarized manifolds with total deficiency one. {I}},
        date={1980},
        ISSN={0025-5645,1881-1167},
     journal={J. Math. Soc. Japan},
      volume={32},
      number={4},
       pages={709\ndash 725},
         url={https://doi.org/10.2969/jmsj/03240709},
      review={\MR{589109}},
}

\bib{Fujita}{article}{
      author={Fujita, Takao},
       title={On polarized varieties of small {$\Delta $}-genera},
        date={1982},
        ISSN={0040-8735,2186-585X},
     journal={Tohoku Math. J. (2)},
      volume={34},
      number={3},
       pages={319\ndash 341},
         url={https://doi.org/10.2748/tmj/1178229197},
      review={\MR{676113}},
}

\bib{Ful98}{book}{
      author={Fulton, William},
       title={Intersection theory},
     edition={Second},
      series={Ergebnisse der Mathematik und ihrer Grenzgebiete. 3. Folge. A Series of Modern Surveys in Mathematics [Results in Mathematics and Related Areas. 3rd Series. A Series of Modern Surveys in Mathematics]},
   publisher={Springer-Verlag, Berlin},
        date={1998},
      volume={2},
        ISBN={3-540-62046-X; 0-387-98549-2},
         url={https://doi.org/10.1007/978-1-4612-1700-8},
      review={\MR{1644323}},
}

\bib{Har77}{book}{
      author={Hartshorne, Robin},
       title={Algebraic geometry},
      series={Graduate Texts in Mathematics, No. 52},
   publisher={Springer-Verlag, New York-Heidelberg},
        date={1977},
        ISBN={0-387-90244-9},
      review={\MR{0463157}},
}

\bib{Hara98}{article}{
      author={Hara, Nobuo},
       title={A characterization of rational singularities in terms of injectivity of {F}robenius maps},
        date={1998},
        ISSN={0002-9327},
     journal={Amer. J. Math.},
      volume={120},
      number={5},
       pages={981\ndash 996},
         url={http://muse.jhu.edu/journals/american_journal_of_mathematics/v120/120.5hara.pdf},
      review={\MR{1646049}},
}

\bib{Hara(two-dim)}{article}{
      author={Hara, Nobuo},
       title={Classification of two-dimensional {$F$}-regular and {$F$}-pure singularities},
        date={1998},
        ISSN={0001-8708},
     journal={Adv. Math.},
      volume={133},
      number={1},
       pages={33\ndash 53},
         url={https://doi.org/10.1006/aima.1997.1682},
      review={\MR{1492785}},
}

\bib{Hara-Watanabe}{article}{
      author={Hara, Nobuo},
      author={Watanabe, Kei-Ichi},
       title={F-regular and {F}-pure rings vs. log terminal and log canonical singularities},
        date={2002},
        ISSN={1056-3911},
     journal={J. Algebraic Geom.},
      volume={11},
      number={2},
       pages={363\ndash 392},
         url={https://doi.org/10.1090/S1056-3911-01-00306-X},
      review={\MR{1874118}},
}

\bib{Kaw2}{article}{
      author={Kawakami, Tatsuro},
       title={On {K}awamata-{V}iehweg type vanishing for three dimensional {M}ori fiber spaces in positive characteristic},
        date={2021},
        ISSN={0002-9947},
     journal={Trans. Amer. Math. Soc.},
      volume={374},
      number={8},
       pages={5697\ndash 5717},
         url={https://doi-org.utokyo.idm.oclc.org/10.1090/tran/8369},
      review={\MR{4293785}},
}

\bib{KT23}{article}{
      author={Kawakami, Tatsuro},
      author={Totaro, Burt},
       title={Endomorphisms of varieties and {B}ott vanishing},
        date={2023},
     journal={preprint available at arXiv:2302.11921v4},
}

\bib{Kawakami-Tanaka(dPsurface)}{article}{
      author={Kawakami, Tatsuro},
      author={Tanaka, Hiromu},
       title={Global ${F}$-regularity for weak del {P}ezzo surfaces},
        date={2024},
     journal={arXiv preprint arXiv:2404.04790},
}

\bib{KT-Fano3}{article}{
      author={Kawakami, Tatsuro},
      author={Tanaka, Hiromu},
       title={Vanishing theorems for {F}ano threefolds in positive characteristic},
        date={2024},
     journal={arXiv preprint arXiv:2404.04764},
}

\bib{KTTWYY1}{article}{
      author={Kawakami, Tatsuro},
      author={Takamatsu, Teppei},
      author={Tanaka, Hiromu},
      author={Witaszek, Jakub},
      author={Yobuko, Fuetaro},
      author={Yoshikawa, Shou},
       title={Quasi-${F}$-splittings in birational geometry},
        date={2022},
     journal={arXiv preprint arXiv:2208.08016v1},
}

\bib{KTY}{article}{
      author={Kawakami, Tatsuro},
      author={Takamatsu, Teppei},
      author={Yoshikawa, Shou},
       title={Fedder type criteria for quasi-${F}$-splitting},
        date={2022},
     journal={arXiv preprint arXiv:2204.10076},
}

\bib{Meg98}{article}{
      author={Megyesi, G.},
       title={Fano threefolds in positive characteristic},
        date={1998},
        ISSN={1056-3911},
     journal={J. Algebraic Geom.},
      volume={7},
      number={2},
       pages={207\ndash 218},
      review={\MR{1620094}},
}

\bib{PW22}{article}{
      author={Patakfalvi, Zsolt},
      author={Waldron, Joe},
       title={Singularities of general fibers and the {LMMP}},
        date={2022},
        ISSN={0002-9327,1080-6377},
     journal={Amer. J. Math.},
      volume={144},
      number={2},
       pages={505\ndash 540},
         url={https://doi.org/10.1353/ajm.2022.0009},
      review={\MR{4401510}},
}

\bib{Spr98}{article}{
      author={Spreafico, Maria~Luisa},
       title={Axiomatic theory for transversality and {B}ertini type theorems},
        date={1998},
        ISSN={0003-889X},
     journal={Arch. Math. (Basel)},
      volume={70},
      number={5},
       pages={407\ndash 424},
         url={https://doi-org.utokyo.idm.oclc.org/10.1007/s000130050213},
      review={\MR{1612610}},
}

\bib{SS10}{article}{
      author={Schwede, Karl},
      author={Smith, Karen~E.},
       title={Globally {$F$}-regular and log {F}ano varieties},
        date={2010},
        ISSN={0001-8708},
     journal={Adv. Math.},
      volume={224},
      number={3},
       pages={863\ndash 894},
         url={https://doi-org.utokyo.idm.oclc.org/10.1016/j.aim.2009.12.020},
      review={\MR{2628797}},
}

\bib{Tan18b}{article}{
      author={Tanaka, Hiromu},
       title={Behavior of canonical divisors under purely inseparable base changes},
        date={2018},
        ISSN={0075-4102},
     journal={J. Reine Angew. Math.},
      volume={744},
       pages={237\ndash 264},
         url={https://doi.org/10.1515/crelle-2015-0111},
      review={\MR{3871445}},
}

\bib{Tan21i}{article}{
      author={Tanaka, Hiromu},
       title={Invariants of algebraic varieties over imperfect fields},
        date={2021},
        ISSN={0040-8735,2186-585X},
     journal={Tohoku Math. J. (2)},
      volume={73},
      number={4},
       pages={471\ndash 538},
         url={https://doi.org/10.2748/tmj.20200611},
      review={\MR{4355058}},
}

\bib{Tan22}{article}{
      author={Tanaka, Hiromu},
       title={Vanishing theorems of {K}odaira type for {W}itt canonical sheaves},
        date={2022},
        ISSN={1022-1824,1420-9020},
     journal={Selecta Math. (N.S.)},
      volume={28},
      number={1},
       pages={Paper No. 12, 50},
         url={https://doi.org/10.1007/s00029-021-00736-0},
      review={\MR{4346509}},
}

\bib{FanoI}{article}{
      author={Tanaka, Hiromu},
       title={Fano threefolds in positive characteristic {I}},
        date={2023},
     journal={arXiv preprint arXiv:2308.08121, to appear in Kyoto J. Math.},
}

\bib{Tan-Bertini}{article}{
      author={Tanaka, Hiromu},
       title={Bertini theorems admitting base changes},
        date={2024},
        ISSN={0021-8693,1090-266X},
     journal={J. Algebra},
      volume={644},
       pages={64\ndash 125},
         url={https://doi.org/10.1016/j.jalgebra.2023.12.038},
      review={\MR{4695615}},
}

\bib{Tot19}{article}{
      author={Totaro, Burt},
       title={The failure of {K}odaira vanishing for {F}ano varieties, and terminal singularities that are not {C}ohen-{M}acaulay},
        date={2019},
     journal={J. Algebraic Geom.},
      volume={28},
      number={4},
       pages={751\ndash 771},
}

\bib{TWY24}{article}{
      author={Tanaka, Hiromu},
      author={Witaszek, Jakub},
      author={Yobuko, Fuetaro},
       title={Quasi-${F^{e}}$-splittings and quasi-${F}$-regularity},
        date={2024},
     journal={arXiv preprint arXiv:2404.06788},
}

\bib{Yob19}{article}{
      author={Yobuko, Fuetaro},
       title={Quasi-{F}robenius splitting and lifting of {C}alabi-{Y}au varieties in characteristic {$p$}},
        date={2019},
        ISSN={0025-5874,1432-1823},
     journal={Math. Z.},
      volume={292},
      number={1-2},
       pages={307\ndash 316},
         url={https://doi.org/10.1007/s00209-018-2198-7},
      review={\MR{3968903}},
}

\bib{Yobuko2}{incollection}{
      author={Yobuko, Fuetaro},
       title={On the {F}robenius-splitting height of varieties in positive characteristic},
        date={2020},
   booktitle={Algebraic number theory and related topics 2016},
      series={RIMS K\^{o}ky\^{u}roku Bessatsu, B77},
   publisher={Res. Inst. Math. Sci. (RIMS), Kyoto},
       pages={159\ndash 175},
      review={\MR{4278982}},
}

\end{biblist}
\end{bibdiv}
